

\documentstyle{amsppt}

\magnification = 1200
\hsize = 6 true in
\baselineskip = 12pt

\redefine\P{\bold P^N}
\define\PP{\bold P^1}
\redefine\L{\Lambda}
\define\St{{\tilde S}}
\define\Xt{{\tilde X}}
\define\Tt{{\tilde T}}
\define\a{\alpha}
\redefine\b{\beta}
\define\BX{\Cal{B}_X}
\define\BXL{\Cal{B}_{X, \L}}

\define\CX{\Cal{C}_X}
\define\CXL{\Cal{C}_{X, \L}}
\define\CXw{\Cal{C}_{X, w}}

\define\tri{\therosteritem}
\define\G{\Gamma}
\define\fwd{\fracwithdelims}
\define\N'{N^{\prime}}
\define\X'{X^{\prime}}
\redefine\H'{H^{\prime}}
\define\CLZ{\text{Cone}_L(Z)}

\define\Lp{L^{\prime}}

\NoRunningHeads
\topmatter
\title Geometric Properties of the Double-Point Divisor \endtitle
\author Bo Ilic \endauthor
\affil Columbia University \endaffil
\endtopmatter

\heading 1. Introduction \endheading

Consider a nonsingular, nondegenerate, projective variety $X^n \subset
\P$ of degree $d$ and codimension $N-n \ge 2$.
One way to study $X$ is to find a new variety $Y$ and a
surjective morphism $f:X \to Y$ and then study $X$ via this
information. In projective geometry, there are two favorite ways of
doing this.

Firstly, we can project $X$ from a generic linear subspace
$\Lambda \subset \bold P^N$ of dimension $N-n-1$ to obtain a finite
surjective morphism $f:X \to \bold P^n$.
This is the higher dimensional analog  of
studying a Riemann surface as a branched cover of the Riemann sphere.
The locus on $X$ where $f$ is not locally an isomorphism is a divisor
$\BXL$ called the {\sl ramification divisor}.
If $H$ is a hyperplane section of $X$ and
$K_X$ is a canonical divisor on $X$
then $\BXL \sim K_X + (n+1)H$. Thus
by varying the generic center of projection $\L$ we
obtain linearly equivalent divisors and so it makes sense to consider
the linear system $|\BX|$ without reference to a center of projection.
L.~Ein proved in \cite{Ein}
that $|\BX|$ is very ample by showing that $|K_X + nH|$ is base point
free. Finally, P. Ionescu, A.J. Sommese, and A. Van de Ven
improved this result:
\proclaim{Adjunction Mapping Theorem} \cite{Io, Theorem 1.4}
Let $X^n \subset \P$ be a nonsingular projective variety. Let
$H$ be a hyperplane section of $X$. If $|K_X + (n-1)H| = \emptyset$
then $X$ is one of the following:
\roster
\item $\P$.
\item The Veronese surface in $\bold P^5$.
\item A quadric hypersurface.
\item The projectivization of a rank $n$ vector bundle on a curve with
the fibers embedded linearly.
\endroster
Moreover, if $|K_X+(n-1)H| \ne \emptyset$ then it is base point free.
\endproclaim

Secondly, we can project $X$ from a generic linear subspace
$\Lambda \subset \bold P^N$ of dimension $N-n-2$ obtaining a finite birational
morphism $f:X \to Y \subset \bold P^{n+1}$ where Y is a hypersurface.
The locus on $X$ where $f$ is not an isomorphism is a divisor
$\CXL$ called the {\sl double-point divisor}. Since
$\CXL \sim (d-n-2)H - K_X$, by varying the generic center of projection
$\L$ we obtain linearly equivalent divisors and so it makes sense to
consider the linear system $|\CX|$ without reference to a center of
projection.  D. Mumford proved that $|\CX|$ is base point free
\cite{BM, Technical Appendix to Section 3, Part 4}. Thus it seems
natural to wonder to what extent a result parallel to Ein's holds in
this case.

For curves and varieties such that $K_X$ is a multiple of a hyperplane
section, such as complete intersections, $|\CX|$ is very ample. However,
in general this does not hold. The main result of this thesis
is a classification of those varieties $X$ such that $|\CX|$ is ample.

To state the classification theorem we need:
\subheading{Definition 3.1} Let $X^n \subset \P$ be a nonsingular,
nondegenerate projective variety of codimension $N-n \ge 2$.
Suppose that there
is an $n+1$ dimensional rational normal scroll
$S = S_{0,0,a_1, \dots, a_{n-1}}$ with all $a_i \ge 1$,
$\sum_{i=1}^{n-1}a_i = N-n$ and singular locus $L$ ($L$ is a line) such
that $L \subset X \subset S$. Then X is a {\sl Roth variety}.

Then the main classification theorem is:
\proclaim{Theorem 4.2} Let $X^n \subset \P$ be a nonsingular, nondegenerate
 variety of codimension $N-n \ge 2$. Let $|\CX|$ denote the
linear system determined by the double-point divisor on $X$.
Then the following are equivalent:
\roster
\item $|\CX|$ is ample.
\item $|\CX|$ separates points.
\item $X$ is not an isomorphic projection of a Roth variety.
\endroster
If, in addition, $X$ is linearly normal then \tri3 can be replaced by:
\roster
\item"$(3^{\prime})$" $X$ is not a Roth variety.
\endroster
\endproclaim

The proof relies on techniques of projective geometry. The main
geometric lemma is a classification of varieties with
certain highly degenerate relative secant varieties:
\proclaim{Theorem 4.3} Let $X^n \subset \P$ be a nonsingular, nondegenerate,
variety of codimension $N-n \ge 2$.
Then the following are equivalent:
\roster
\item There exists a positive dimensional subvariety
$W$ of $X$ such that the secant variety $S(W,X)$ has dimension $n+1$.
\item There is $p \ne q \in X$ such that $S(p,X) = S(q,X)$.
\item $X$ is an isomorphic projection of a Roth variety.
\endroster
If, in addition, $X$ is linearly normal then \tri3 can be replaced by:
\roster
\item"$(3^{\prime})$" $X$ is a Roth variety.
\endroster
\endproclaim
See the beginning of Section 2 for the definition of the relative
secant varieties.

We provide a convenient characterization and existence
theorem for Roth varieties:
\proclaim{Theorem 3.8}
Let $S^{n+1} \subset \P$ be a rational normal scroll
$S_{0,0,a_1,\dots,a_{n-1}}$ with all $a_i \ge 1$ and $\sum_{i=1}^{n-1}a_i =
N-n$ and let $L$ be the singular locus of $S$. Let
$\Cal E = \Cal O_{\PP} \oplus \Cal O_{\PP} \oplus \Cal
O_{\PP}(a_1) \oplus \dots \oplus \Cal O_{\PP}(a_{n-1})$ and let
$\pi_1:\bold P(\Cal E^*) \to \PP$ be a $\bold P^n$ bundle with diagram:
$$
\CD
\bold P(\Cal E^*) @>\pi_2>> S \subset \P \\
@V\pi_1VV @. \\
\PP @.
\endCD
$$
that desingularizes $S$ such that the morphism $\pi_2$ is given by the
complete linear series $|\Cal O_{\bold P(\Cal E^*)}(1)|$. Let H be the
pullback of a hyperplane in $\P$ via $\pi_2$ and let $F$ be a fiber of $\pi_1$.
Then for every $b > 0$ there exists a nonsingular variety $\Xt \in |bH +F|$
such that if $X = \pi_2(\Xt)$ then $\pi_2|_{\Xt}:\Xt \to X$ is an
isomorphism and $X$ is a linearly normal Roth variety of
degree $b(N-n)+1$ (with associated scroll $S$).
Conversely, if $X$ is a Roth variety of degree $d$ such that
$L \subset X \subset S$, then $b = \frac{d-1}{N-n}$ is an integer and there
is a desingularization of $S$ as above; $\pi_2: \bold P(\Cal E^*) \to S$ and
$\Xt \in |bH + F|$ which is mapped isomorphically onto $X$ by $\pi_2$.
(In particular, Roth varieties are linearly normal).
\endproclaim

We determine further properties of Roth varieties and in
particular their relationship to Castelnuovo varieties. Finally,
an application to vanishing related to the regularity conjecture is made.

\subheading{Historical background}

Roth surfaces were studied by L. Roth in \cite{Ro, \S 3.5}
where he proved, for
example, that such surfaces have maximal sectional genus.  They were
re-examined in modern language by A.J. Sommese in \cite{So, \S 1.1}
in order to
settle ``the main technical problem'' in his first paper on the
adjunction mapping. Sommese computed the sectional genus of a Roth
surface and gave a geometric description of Roth surfaces as in
Proposition 5.1 \tri2.

Restated in the language of projective geometry, he also
proved:
\proclaim{Lemma 1.1.1} Let $X \subset \P$ be a nonsingular,
linearly normal surface. Let $L \subset X$ be a line and suppose that
$H^1(X, \Cal O_X(K_X + H - L)) \ne 0$. Then projection from the line
$L$ extends to a morphism $\pi_L : X \to \bold P^{N-2}$ whose image is
a rational normal curve of degree $N-2$.
\endproclaim
{}From this
theorem it is easy to deduce that $X$ is a Roth surface.  This is the
first theorem where the internal properties of $X$ were used
to deduce the existence of a rational normal scroll $S_{0,0,N-2}$
containing $X$.  Proposition 4.11 of this work follows the main line of
argument of Sommese's proof of the above lemma.

Roth surfaces contained in $\bold P^4$ were extensively studied by
K. Hulek, C. Okonek, and A. Van de Ven in \cite{HOV}.
Let $X \subset \bold P^4$ be a smooth surface. An ideal subsheaf
$\Cal I_{\Xt} \subset \Cal I_{X}$ determines a {\sl multiplicity-2 structure}
on $X$ if for each $p \in X$ there are local coordinates $x_0, \dots, x_3$
such
that $\Cal (I_X)_p = (x_0, x_1)$ and $\Cal (I_{\Xt})_p = (x_0,x_1^2)$.
Their main theorem, somewhat restated, is:
\proclaim{Theorem 13}
Let $X^2 \subset \bold P^4$ be a nonsingular Castelnuovo surface of
degree $2b+1$. Then $X$ is a Roth surface $\iff$ there is a multiplicity-2
structure $\Xt$ on $X$ such that $\Xt$ is a complete intersection of
type $(2, 2b+1)$.
\endproclaim
They also proved that for every positive integer $b$ there exists a
Roth surface of degree $2b+1$ contained in $\bold P^4$ (Proposition
10). The proof of the general existence result in this paper (Theorem 3.8)
is not a generalization of their proof.

\subheading{Acknowledgements} I thank Professor Henry Pinkham for his
help and encouragement throughout my graduate studies at Columbia
University where this work was prepared. I also thank Professor Robert
Lazarsfeld for discussing an early version of this paper with me.

\newpage

\subheading{Conventions and basic terminology}

1. The ground field is $\bold C$, the field of complex numbers.

2. Unless otherwise stated, varieties are assumed to be irreducible and
projective.

3. Denoting a variety by $X^n$ means that $X$ is a variety of dimension $n$.

4. A variety $X \subset \P$ is {\sl nondegenerate} if it is not contained in
any hyperplane $H \subset \P$.

5. If $X$ and $Y$ are varieties such that $Y \subset X$ then
$\Cal I_{Y/X}$ denotes the sheaf of ideals of $Y$ in $X$. We
write $\Cal I_Y$ for $\Cal I_{Y/\P}$.

6. If $X^n \subset \P$ is a nondegenerate variety then $X$ is {\sl linearly
normal} if one of the following equivalent conditions holds:
\roster
\item There does not exist a nondegenerate variety $\X' \subset
\bold P^{N+1}$ and a point $p \in \bold P^{N+1} \setminus \X'$ such
that projection from the point $p$, $\pi_p:\bold P^{N+1} - - \to \P$
restricts to an isomorphism $\pi_p|_{\X'}: \X' \to X$.
\item If $H$ is the restriction of a hyperplane in $\P$ to $X$ then
$X \subset \P$ is an embedding of $X$ given by the complete
linear system $|H|$.
\item $H^1(\P, \Cal I_X(1)) = 0.$
\item The natural map $H^0(\P, O_{\P}(1)) \to H^0(X, \Cal O_X(1))$
is surjective.
\endroster
$X$ is {\sl k normal} if $H^1(\P, \Cal I_X(k)) = 0$.
$X$ is {\sl k regular} if $H^i(\P, \Cal I_X(k-i)) = 0$ for all $i >0$.
Notice that if $X$ is $k$ regular then $X$ is also $k-1$ normal.
$X$ is {\sl projectively normal} if $X$ is normal and $k$ normal for
all $k \ge 0$. $X$ is {\sl arithmetically Cohen-Macaulay} if
$H^i(\P, \Cal I_X(k)) = 0$ for $1 \le i \le n$ and all $k \in \bold Z$.

7. If $X$ is a variety then $A(X)$ denotes the Chow ring of $X$ i.e.
cycles on $X$ modulo rational equivalence. $A(X) = \oplus_{i=0}^{n}A_i(X)$
where $A_i(X)$ is the dimension $i$ cycles modulo rational equivalence.
If $x \in A_0(X)$ is a
zero cycle then the degree of $X$ is deg($x$) or $\int x$. However, we
usually abuse notation and use $x$ to denote both the zero cycle
and its degree depending on context. If
$x,y \in A(X)$ then $x \cdot y$ is the intersection product of $x$ and $y$.
If $x$ and $y$ are cycles then $x \sim y$ if $x = y$ in $A(X)$.
Finally, $x^k = x \cdot x \cdot \dots \cdot x$ ($k$ times).

8. If $V$ is a vector space then $\bold P(V)$ is the projective space
of lines in $V$. The same convention extends to the vector bundle $\Cal V$ and
its associated projective bundle $\bold P(\Cal V)$.

9. $K_X$ denotes a canonical divisor on the variety $X$. If $X \subset \P$
then $H$ denotes a hyperplane restricted to $X$.

10. If $X^n \subset \P$ is a variety and $p \in X$ is a nonsingular
point of $X$ then $T_pX$ is the Zariski tangent space to $X$ at $p$ and
$\Tt_pX$ is the embedded projective tangent space to $X$ at $p$.
(Hence $T_pX$ is an $n$ dimensional complex vector space and $\Tt_pX$ is
an $n$ dimensional projective linear subspace of $\P$).

11. If $X$ is a subset of $\P$ then $<X>$ is the smallest linear
subspace of $\P$ that contains $X$. If $X$ and $Y$ are subsets of $\P$ then
$<X,Y> = <X \cup Y>$.

12. The notation for rational normal scrolls and their basic projective
geometric properties are taken from \cite{Ha2}. Intersection theory
on a rational normal scroll $S \subset \P$ and the description of $S$
as the image of a morphism $\St \to S$ where $\St$ is a
projective bundle over $\PP$ are
discussed in \cite{Ha1} and \cite{EH}. For cohomological calculations on
scrolls, see \cite{Hart1, Exercises III.8.1, III.8.4}.

\newpage

\heading 2. Double-point divisors \endheading

\subheading{Definition and basic properties of relative secant varieties}

Let $X^n \subset \P$ be a variety and let $W^m \subset X$ be a subvariety.
Following F. Zak's terminology, we define the {\sl relative secant variety}
$S(W,X)$ as the special case of the join of $W$ and $X$ when
$W \subset X$. Our definition is the specialization of the definition
of join found in \cite{Ha2, Example 8.1}. Let
$\Delta = \{ (p,q) \in W \times X \mid p = q \}$ be the diagonal.
Consider the morphism
$$
\matrix
W \times X \setminus \Delta & @>\phi>> & \bold G(1, \P) \\
(p,q) & & <p,q>
\endmatrix
$$
and let $\Cal S(W,X) = \overline{\phi(W\times X \setminus \Delta)}$.
Let $I \subset \P \times \bold G(1, \P)$ be the incidence correspondence
$$
I = \{ (p,l) \mid p \in l, \  l \in \Cal S(W,X) \}
$$
with natural projections
$$
\CD
I @>\phi_2>> \Cal S(W,X) \\
@V\phi_1VV @. \\
\P @.
\endCD
$$
then $S(W,X) = \phi_1(I)$ is an irreducible variety and
$\text{dim}(S(W,X)) \le m + n + 1$. Also, $\text{dim}(S(W,X)) \ge n$ and
$\text{dim}(S(W,X)) = n$ if and only if $X$ is a linear space.
$S(X,X)$ is the usual secant variety of $X$. If $W$ is a possibly
reducible subvariety of $X$ with irreducible decomposition
$W = W_1 \cup \dots \cup W_k$ then
$S(W,X) = S(W_1, X) \cup \dots \cup S(W_k,X)$. If $W_1 \subset W_2$ are
possibly reducible subvarieties of $X$ then $S(W_1,X) \subset S(W_2,X)$.

If $p \in X$ we write $S(p,X)$ for $S(\{p\},X)$. If $X$ is not a linear
space then $\text{dim}(S(p,X)) = n+1$. If $p$ is a nonsingular point
of $X$ then
$S(p,X) = \{ r \in \P \mid r$ \ lies on a proper bisecant of $X$ through
$p$ i.e. there exists $q \in X, q \ne p$ \ such that $r \in <p,q> \}
\cup \Tt_pX$. In particular, $\Tt_pX \subset S(p,X)$.
Moreover, if $H$ is a hyperplane that meets $X$ transversally at a
nonsingular point $p$ then $S(p, X \cap H) = S(p, X) \cap H$.
If $p \in X$ then by repeatedly taking
hyperplane sections through $p$ we see that $\text{deg}(S(p,X)) <
\text{deg}(X)$.

\subheading{Basic properties of the double-point divisor}

Let $X^n \subset \P$ be a nonsingular, nondegenerate variety of degree
$d$ and codimension $N-n \ge 2$. Let $\L$ be a linear subspace of
$\P$ of dimension $N-n-2$ such that $X \cap \L = \emptyset$.  Let
$\pi_{\L}: X \to \bold P^{n+1}$ denote the morphism obtained by restricting
the projection with center $\L$, $\P - - \to \bold P^{n+1}$ to $X$.
$\pi_{\L}$ is a finite morphism and so $Y = \pi_{\L}(X)$ is a hypersurface
in $\bold P^{n+1}$. $\pi_{\L}$ need not be birational onto its image but
for a generic choice of $\L$ it is.

If $p \in X$ and $q = \pi_{\L}(p)$ then we say that $\pi_{\L}:X \to Y$ is
an {\sl isomorphism at} $p$ if there is an open neighborhood $U$ of $q$
such that
$$
\pi_{\L}|_{\pi_{\L}^{-1}(U)} : \pi_{\L}^{-1}(U) \to U
$$
is an isomorphism.
$\pi_{\L}$ is an isomorphism at $p$ if and only if $<\L, p> \cap X = \{p\}$
and $\Tt_pX \cap \L = \emptyset$. The first condition says that
$\pi_{\L}^{-1}(q)$ is a single point and the second one says that
$\pi_{\L}$ does not ramify at $p$. Equivalently, $\pi_{\L}$ is an
isomorphism at $p$ if and only if $S(p,X) \cap \L = \emptyset$.

If $\pi_{\L}$ is birational there is a divisor $\CXL$ on $X$, the {\sl
double-point divisor}, which depends on $X$ and $\L$ such that:
\roster
\item
$$
\align
\text{support}(\CXL) &= \{ p\in X \mid \pi_{\L} \ \text{is not an
isomorphism at \ } p\} \\
&= \{ p \in X \mid S(p,X) \cap \L \ne \emptyset \}
\endalign
$$
\item $\CXL \sim (d-n-2)H - K_{X}$
\item The sheaf of $\Cal O_Y$ ideals $(\pi_{L})_{*}\Cal O_X(-\CXL)$ is
the {\sl conductor} of the finite birational map $\pi_{\L} : X \to Y$.
\endroster
For these basic properties of the double-point divisor see \cite{BM,
Technical Appendix to Section 3, Part 4}
and \cite{Kl, Sections I.D and V.A}.

By \tri2 the divisor class of $\CXL$ does not depend on the center of
projection $\L$ chosen so it makes sense to consider the linear system
$|\CX|$ without reference to a given $\L$.

\subheading{Remark 2.1} If $X^n \subset \P$ is a  variety, $\L^m$ is a center
of projection such that $\L \cap X = \emptyset$ and
$\pi_{\L}: X \to \bold P^{N-m-1}$ is the projection then let $Y = \pi_{\L}(X)$.
If $\pi_{\L}:X \to Y$ is an isomorphism we say that $Y$ is an
{\sl isomorphic projection} of $X$. If this is so then $\Cal O_X(\CX) \cong
\Cal O_Y(\Cal C_{Y})$. Thus, for example,  $X^n \subset \P$ is
a nonsingular, nondegenerate variety such that $|\CX|$ separates points
$\iff$ $X$ is an isomorphic projection of a linearly normal
variety $\Xt$ such that $|\Cal C_{\Xt}|$ separates points.

\proclaim{Proposition 2.2} \cite{BM, Technical Appendix to Section 3, Part 4}
 Let $X^n \subset \P$ be a nonsingular,
nondegenerate variety of codimension $N-n \ge 2$. Then $|\CX|$ is
base point free.
\endproclaim
\demo{Proof}
Given $p \in X$, we can pick an $N-n-2$ plane $\L$ such that
$\L \cap S(p,X) = \emptyset$. Then if $\pi_{\L}: X \to \bold P^{n+1}$ is
projection from the center $\L$ and $Y = \pi_{\L}(X)$ then
$\pi_{\L}:X \to Y$ is a finite, birational morphism and is an isomorphism
at $p$. Hence $\CXL$ is defined and $p \notin \CXL$. Thus, $p$ is not
a base point of $|\CX|$.
\qed
\enddemo

\proclaim{Corollary 2.3}
$|\CX|$ is nef i.e. $\CX \cdot W \ge 0 $ for all irreducible curves
$W \subset X$.
\endproclaim
\demo{Proof}
Let $W \subset X$ be an irreducible curve and let $p \in W$. As above,
pick $\L$ such that $\CXL$ is defined and $p \notin \CXL$. Then
$\CXL \cdot W \ge 0$ since $\CXL$ is effective.
\qed
\enddemo

\proclaim{Lemma 2.4}
If $p,q \in X$ and $S(p,X) \ne S(q,X)$ then there exists an $N-n-2$ plane
$\L$ such that $\L \cap S(p,X) = \emptyset$ but $\L \cap S(q,X) \ne
\emptyset$. (So that $\CXL$ is defined and $p \notin \CXL$ but
$q \in \CXL$).
\endproclaim

\demo{Proof}
Let $M$ be an $N-n-1$ plane that intersects both $S(q,X)$ and $S(p,X)$ in a
finite collection of points (the generic $N-n-1$ plane does this) and
such that there is a point $r \in (M \cap S(q,X)) \setminus S(p,X)$.
If we cannot do this then $S(q,X) \subset S(p,X)$ which is a contradiction.
Let $\L$ be the $N-n-2$ plane determined by a hyperplane in $M$ which
contains $r$ but avoids the finite set of points $M \cap S(p,X)$. $\L$
satisfies the requirements of the lemma.
\qed
\enddemo

\proclaim{Corollary 2.5}
 If $|\CX|$ does not separate the points $p$ and $q$ then
$S(p,X) = S(q,X)$.
\endproclaim

To prove the next result we need to use the Kleiman-Nakai-Moishezon
criterion for ampleness:
\proclaim{Theorem 2.6} \cite{Hart2, Chapter I, Theorem 5.1}
Let $X$ be a nonsingular projective variety and let $L$ be a line
bundle on $X$. Then $L$ is ample $\iff L^{\text{dim}(W)} \cdot W > 0$
for every subvariety $W$ of $X$ with $\text{dim}(W) \ge 1$.
\endproclaim

\proclaim{Corollary 2.7}
If $|\CX|$ separates points then $|\CX|$ is ample.
\endproclaim

We'll use Zak's theorem on tangencies in proving the next
two propositions:
\proclaim{Theorem 2.8} \cite{Zak, Chapter 1, Corollary 1.8}
If $X^n \subset \P$ is a nonsingular, nondegenerate variety and
$L \subset \P$ is an $m$ dimensional linear subspace with
$n \le m \le N-1$ then
$\text{dim} \{p \in X \mid \Tt_pX \subset L \} \le m-n.$
\endproclaim

\proclaim{Proposition 2.9} If $W$ is a positive dimensional
subvariety of $X$ such that $S(p,X) \allowmathbreak= S(W,X)$ for all  $p \in W$
then
$S(W,X) = \cup_{p \in W} \Tt_pX$.
\endproclaim
\demo{Proof}
$\Tt_pX \subset S(W,X)$ for all $p \in W$. So $T = \cup_{p \in W}\Tt_pX$ is
subvariety of $S(W,X)$. Since $S(W,X)$ is $n+1$ dimensional and
irreducible, if $T \ne S(W,X)$ then $\Tt_pX = \Tt_qX$ for all $p,q \in W$.
Since $\text{dim}(W) \ge 1$ this contradicts Zak's theorem.
\qed
\enddemo

\proclaim{Proposition 2.10} If $W$ is a possibly
reducible subvariety of $X$ such that $S(p,X) = S(W,X)$ for all
$p \in W$ then $W \subset  \cap_{p \in W} \Tt_pX$.
Moreover, if $W$ is positive dimensional then
 $\text{dim}(\cap_{p \in W} \Tt_pX) \le n-1$.
\endproclaim

\demo{Proof}
Pick $p \in W$. We'll show that $W \subset \Tt_pX$. Let $q \in W$.
Suppose that $q \notin \Tt_pX$. Then we can ``define'' a map:
$$
\Tt_pX \setminus (\Tt_qX \cup X) \to <\Tt_pX,q>
$$
where a point $r$ goes to a point $w \in X$, $w \ne q$ on the line
$<q,r>$. Observe that:
\roster
\item For each $r \in \Tt_pX$, $r \in S(p,X)$ and thus $r \in S(q,X)$ and
since $r \notin \Tt_qX$, $<q,r>$ is a proper bisecant of $X$.
\item The map is not well defined since there could be more than one
choice of $w$, but complex analytic locally, the map can be defined
to be analytic; use the implicit function theorem on the natural
incidence relation.
\item The map is injective.
\endroster

Thus we have an $n$ dimensional complex neighborhood of $X \subset
<\Tt_pX,q>$. Thus, $X \subset <\Tt_pX,q>$ which contradicts the
nondegeneracy of $X$. Thus, $q \in \Tt_pX$. Since $q \in W$ was
arbitrary, $W \subset \Tt_pX$. Since $p \in W$ was arbitrary the
first part of the proposition follows.

If $\text{dim}(\cap_{p \in W} \Tt_pX) = n$ then
$\Tt_pX = \Tt_qX$ for all $p,q \in W$. But since $\text{dim}(W) \ge 1$
this contradicts Zak's theorem.
\qed
\enddemo

\proclaim{Proposition 2.11}
Let $X \subset \P$ be a nondegenerate, nonsingular curve of
codimension $N-1 \ge 2$. Then $|\CX|$ is very ample.
\endproclaim

\demo{Proof}
By taking a generic projection, we can assume without loss of generality
that $X \subset \bold P^3$ (see Remark 2.1).
First we'll show that $|\CX|$ separates
points. Pick $p,q \in X$, $p \ne q$. Pick $r \in X$ not on the line
$<p,q>$ and consider the line $L = <q,r>$. If there exists $w \in
L \setminus X$ such that $w \notin S(p,X)$ then $\CXw$ is  defined
and $p \notin \CXw$ but $q \in \CXw$.  Thus suppose that $w \in S(p,X)$
for all $w \in L \setminus X$. Then for generic $w \in L$, $<w,p>$ is
a proper bisecant of $X$. Thus varying $w$ we see that the plane
$<L,p>$ contains infinitely many points of $X$ so that
$X \subset <L,p>$. But this contradicts the nondegeneracy of $X$.

Now we'll show that $|\CX|$ separates tangent directions i.e. given
$q \in X$ there exists $w \in \P$ such that $|\CXw|$ is defined,
$q \in \CXw$ but  $w \notin \Tt_qX$. Clearly we
can pick $r \in X$ such that $<q,r>$ is a bisecant of $X$ and $r
\notin \Tt_qX$. Then there must be some $w \in <q,r>$ that works or
we get a contradiction as above.
\qed
\enddemo

The last proposition can also be proven by using Castelnuovo's bound for
the genus of a space curve \cite{Hart1, Theorem IV.6.4}.
If $g$ is the genus of $X$ then
one shows that $\text{deg}(\CX) = \text{deg}((d-3)H- K_X) \ge 2g+1$ and
so $\CX$ is very ample.

Recall that a variety $X^n \subset \P$ is {\sl semi-canonical} if $K_X \sim
\alpha H$ for some integer $\alpha$. For example, by the Barth-Larsen theorem
\cite{La}, if $n \ge (N-n) + 2$ then
$X$ is semi-canonical. Also, by using the adjunction formula,
we see that complete intersections are semi-canonical.
\proclaim{Proposition 2.12}
Let $X^n \subset \P$ be a nonsingular, nondegenerate, semi-canonical variety of
codimension $N-n \ge 2$. Then $|\CX|$ is very ample.
\endproclaim
\demo{Proof}
Let $d$ be the degree of $X$. Then $\CX \sim (d-n-2)H-K_X \sim \alpha H$ for
some integer $\alpha$. We can choose an $N-n-2$ plane $\L$ in $\P$ such
that $\CXL$ is an element of $|\CX|$. Hence $H^0(X, \Cal O_X(\CX)) \ne 0$
and so $\alpha \ge 1$. Thus $|\CX|$ is very ample.
\qed
\enddemo

\newpage

\heading 3. Roth varieties \endheading

\subheading{Definition 3.1} Let $X^n \subset \P$ be a nonsingular,
nondegenerate projective variety of codimension  $N - n \ge 2$.
Suppose that there
is an $n+1$ dimensional rational normal scroll
$S = S_{0,0,a_1, \dots, a_{n-1}}$ with all $a_i \ge 1$,
$\sum_{i=1}^{n-1}a_i = N-n$ and singular locus $L$ ($L$ is a line) such
that $L \subset X \subset S$. Then X is a {\sl Roth variety}.

Let $S$ be the rational normal scroll $S_{0,0,a_1, \dots, a_{n-1}}
\subset \P = \bold P(V)$ where $V$ is an $N+1$ dimensional vector space
and all $a_i \ge 1$. Then
$\text{dim}(S) = n +1$, $d_S := \text{deg}(S) = \sum_{i=1}^{n-1} a_i$, and
$N = d_S + n$. Let $L$ be the singular locus of $S$. $L$ is a line contained
in $S$. Then $S$ can be considered as the cone with vertex $L$ over the
nonsingular rational normal scroll $S_{a_1, a_2, \dots, a_{n-1}}$
where $<S_{a_1, a_2, \dots, a_{n-1}}>$ and $L$ are
complementary linear subspaces of $\P$.

$S$ can also be described geometrically as follows:
there are linear subspaces $V_0, V_1, \dots, V_{n-1}$ of
$V$ such that:
\roster
\item $V_i \cap V_j = 0$ for all $i$ and $j$
\item $V = V_0 \oplus V_1 \oplus \dots \oplus V_{n-1}$
\item $\text{dim} (V_0) = 2$
\item $\text{dim} (V_i) = a_i+1$ if $i \ge 1$
\endroster
and such that for each $i \ge 1$ we can find a rational normal curve
$C_i$ of degree $a_i$ in $\bold P(V_i)$ and isomorphisms $\phi_i: \PP
\to C_i$ such that:
$$S = \bigcup_{p \in \PP} <L, \phi_1(p), \phi_2(p),
\dots , \phi_{n-1}(p)>.$$
Note also that if $p \ne q$ then: $$<L,
\phi_1(p), \phi_2(p), \dots , \phi_{n-1}(p)> \cap <L, \phi_1(q),
\phi_2(q), \dots , \phi_{n-1}(q)> = L.$$

Let $\St \subset \PP \times \P$ be defined by:
$$\St = \{(p, <L, \phi_1(p), \phi_2(p), \dots \phi_{n-1}(p)>)
\mid p \in \PP \}.$$
$\St$ is a variety and if
$$
\CD
\St @>\pi_2>> \P \\
@V\pi_1VV   @.   \\
\PP @. {}
\endCD
$$
are the natural projections then $\pi_1: \St \to \PP$ makes $\St$
into a $\bold P^n$-bundle over $\PP$ and hence $\St$ is
nonsingular. Note also that $\pi_2(\St) = S$. Moreover, $\pi_2: \St
\to S$ is such that $\St \setminus \pi_2^{-1}(L) \to S \setminus L$ is
an isomorphism. Thus $\pi_2$ is birational onto its image. Finally,
$\pi_2^{-1}(L) = \PP \times L \subset \St$.

Let $F$ denote a fiber of $\pi_1$ and $H$ denote the pullback of a hyperplane
section given by the map $\pi_2: \St \to \P$.
In fact, if $\Cal E = \Cal O_{\PP} \oplus \Cal O_{\PP} \oplus
\Cal O_{\PP}(a_1) \oplus \dots \oplus \Cal O_{\PP}(a_{n-1})$ then
$\St = \bold P(\Cal E^*)$,
$\pi_1$ is the natural projection to $\PP$, and $\Cal O_{\bold P(\Cal E^*)}(1)
= \Cal O_{\St}(H)$. Let $A(\St)$ be the
ring of cycles mod rationally equivalent cycles (the Chow ring). Then
$$A(\St) = \bold Z[F,H]/(F^2, H^{n+1}-d_S H^n \cdot F, H^{n+2},
H^{n+1} \cdot F).$$
Also $\text{deg}(H^n \cdot F) = 1$ and $\text{deg}(H^{n+1}) = d_S$.
Recall also that $$K_{\St} \sim -(n+1)H + (d_S - 2)F.$$
\cite{EH} is a reference for these facts.

\proclaim{Lemma 3.2}
In $A(\St)$,  $\PP \times L \sim H^{n-1} - d_S H^{n-2} \cdot F.$
\endproclaim

\demo{Proof}
We know that $\PP \times L \sim \a H^{n-1} + \b H^{n-2} \cdot F$ for some
integers $\a$ and $\b$. Now by the construction we see that
$(\PP \times L) \cdot F \cdot H = 1$. Thus $1 = \a H^n \cdot F + \b H^{n-1}
\cdot F^2 = \a H^n \cdot F = \a$ so that $\a = 1$. Also $(\PP \times L)
\cdot H^2 = 0$ since $\pi_2(\PP \times L) = L$.
Thus, $0 = H^{n+1} + \b H^n \cdot F$ and so $\b = -d_S$.
\qed
\enddemo

We have the diagram of projections obtained by restriction to $\PP \times L$
$$
\CD
\PP \times L @>\pi_2|_{\PP \times L}>> L \\
@V\pi_1|_{\PP \times L}VV @. \\
\PP @. {}
\endCD
$$
The fibers over $\PP$ are linearly equivalent in $A(\PP \times L)$ and thus are
rationally equivalent in $A(\St)$. The same holds for fibers over $L$. Let
$B$ be a fiber over $\PP$ and let $C$ be a fiber over $L$.

\proclaim{Lemma 3.3}
In $A(\St)$, $B \sim H^{n-1} \cdot F$ and
$C \sim H^n - d_S H^{n-1} \cdot F$.
\endproclaim

\demo{Proof}
$B \sim (\PP \times L) \cdot F = (H^{n-1} - d_S H^{n-2} \cdot F) \cdot F =
H^{n-1} \cdot F$ and  $C \sim (\PP \times L) \cdot H
 = (H^{n-1} - d_S H^{n-2}
\cdot F) \cdot H$.
\qed
\enddemo

Let $b$ be a positive integer. Consider the linear system $|bH+F|$ on
$\St$. $|H|$ is clearly base point free. $|H+F|$ is very ample since
$H+F$ is the pullback of a hyperplane section of the projective embedding
of $\St$ as the rational normal scroll $S_{1,1,a_1+1,\dots,a_{n-1}+1}$.
Thus $|bH+F|$ is very ample. Hence, by Bertini's theorem, a generic
$\Xt \in |bH+F|$ is a smooth and irreducible $n$ dimensional variety
and intersects $\PP \times L$ in an irreducible curve.

Now choose an $\Xt$ as above and let $X = \pi_2(\Xt)$.

\proclaim{Proposition 3.4}
The restriction of $\pi_2$ to $\Xt$ gives an
isomorphism $\Xt \cong X$.
\endproclaim

\demo{Proof}
The restriction of $\pi_2$ gives an isomorphism $\Xt \setminus (\PP \times L)
\to X \setminus L$. I claim that $\pi_2|_{\Xt}: \Xt \to X$ is a bijection:
If not, there is a point $p \in L$ such that $(\PP \times p) \cap \Xt$
contains more than one point. But,
$$\Xt \cdot C = (bH+F) \cdot (H^n
- d_S H^{n-1} \cdot F) =bH^{n+1} - b d_S H^{n}\cdot F + H^n \cdot F =
1.$$
So either $\PP \times p$ intersects $\Xt$ in exactly one point $q$
and the intersection is transversal at $q$ or $\PP \times p$ is entirely
contained within $\Xt$.  Suppose the second case occurs. Then $(\PP
\times L) \cdot \Xt = \PP \times p \sim C$ since we know that $\PP \times
L$ meets $\Xt$ in an irreducible curve. Now, $$(\PP \times L) \cdot \Xt
\cdot H = (H^{n-1} - d_S H^{n-2} \cdot F) \cdot (bH+F) \cdot H =
bH^{n+1} + (1-bd_S)H^{n} \cdot F = 1.$$ On the other hand, $C \cdot
H = 0.$ This contradiction shows
that the second case doesn't occur and so our map is a bijection.

Take any point $p \in L$ and let $q$ be the unique point of
intersection of $\PP \times p$ with $\Xt$. To show that our map is an
isomorphism it suffices to show that $\pi_2|_{\Xt}$ is immersive at $q$
i.e. that $d\pi_2|_{\Xt}$ is injective at $q$: In fact, the kernel of
$(d\pi_2) |_q : T_q\St \to T_{\pi_2(q)}\P$ is exactly the
one-dimensional subspace $T_q(\PP \times p)$: clearly
$T_q(\PP \times p)$ is contained in the kernel since $\pi_2$ contracts $\PP
\times p$. On the other hand, $\pi_2$ acts as the identity map on the
$\bold P^n = \pi_1^{-1}(\pi_1(q))$. Thus, since $\PP \times p$ meets $\Xt$
transversally at $q$ our claim follows.
\qed
\enddemo

\proclaim{Proposition 3.5} $X$ is a nonsingular, nondegenerate, and linearly
normal variety.
\endproclaim

\demo{Proof}
$X$ is nonsingular since it is isomorphic to $\Xt$ and $\Xt$ is nonsingular.
A hyperplane containing $X$ pulls back to give a nonzero global section of
$\Cal I_{\Xt / \St}(H)$ where $\Cal I_{\Xt / \St}$ denotes the sheaf of
ideals of $\Xt$ in $\St$. Thus, to prove that $X$ is nondegenerate it
suffices to show that $H^0(\St, \Cal I_{\Xt / \St}(H)) = 0$. Assuming that $X$
is nondegenerate, to show that $X$ is linearly normal it suffices to show that
the morphism $\Xt \to X \subset \P = \bold P(H^0(\St, \Cal O_{\St}(H)))$ is
given by a complete linear series. This follows if
$H^0(\Xt, \Cal O_{\Xt}(H)) = H^0(\St, \Cal O_{\St}(H))$. From the long
exact sequence obtained from the short exact sequence
$$ 0 \to \Cal I_{\Xt / \St}(H) \to \Cal O_{\St}(H) \to \Cal O_{\Xt}(H) \to 0$$
this follows if $H^1(\St, \Cal I_{\Xt / \St}(H)) = 0$. Since
$\Cal I_{\Xt / \St} = \Cal O_{\St}(-bH-F)$ it suffices to show that
$H^i(\St, \Cal O_{\St}((1-b)H - F)) = 0$ for $i = 0,1$ and $b \ge 1$.

This divides into two cases: $b \ge 2$ and $b=1$. So suppose that $b \ge 2$.
By Serre duality it suffices to show that
$H^{i}(\St , \Cal O_{\St}(K_{\St} + ((b-1)H+F))) = 0$ for $i = n, n+1$. Since
$b \ge 2$, $|(b-1)H+F|$ is very ample and so the required vanishing follows
from Kodaira's vanishing theorem.

If $b=1$ the above argument doesn't apply, but a direct cohomology
computation for scrolls does the trick. I want
to show that $H^i(\St, \Cal O_{\St}(-F)) = 0$ for $i = 0,1$.
But $\pi_{1 *} \Cal O_{\St}(-F) = \Cal O_{\PP}(-1)$ and $R^i\pi_{1 *}
\Cal O_{\St}(-F) = 0$ for $i > 0$ and so
$H^i(\St, \Cal O_{\St}(-F)) = H^i(\PP, \Cal O_{\PP}(-1)) = 0$ for $i = 0,1$
as required.
\qed
\enddemo

Let $X^n \subset \P = \bold P(V)$ be a nonsingular, nondegenerate projective
variety.
Let $L^m = \bold P(V_0)$ be a linear space contained in $X$ and let
$\pi_L: X - - \to \bold P(V / V_0)$ be the rational map given by projection
from the center $L$. Let $\G \subset X \times \bold P(V/V_0)$ be the
graph of $\pi_L$. We obtain the diagram:
$$
\matrix
E & \subset & \G & @>f_2>> & \bold P(V/V_0) \\
\Big\downarrow &  & f_1\Big\downarrow \  &  & \Big\Vert \\
L & \subset & X & -- \to & \bold P(V/V_0)
\endmatrix
$$
where $\G$ is the blowup of $X$ along the center $L \subset X$ and
$E$ is the exceptional divisor of the blowup. See [Ha2, Example 7.18].
Alternatively, observe that the graph of the projection from
$\L$, $\pi: \P -- \to \bold P(V/V_0)$ is the blowup of $\P$ along $\L$.
This holds since if $\P$ has coordinates $x_0, \dots, x_N$,
$\L$ is given by the equations $x_{m+1} =  \dots = x_N = 0$ and
$\bold P(V/V_0)$ has coordinates $y_{m+1}, \dots, y_N$ then both
the graph and the blowup are
$$
\{ ((x_0, \dots, x_N), (y_{m+1}, \dots, y_N)) \mid x_i y_j = x_j y_i
   , \ m+1 \le i, j \le N \} \subset \P \times \bold P(V/V_0).
$$
Then the graph of $\pi|_L$ is just the proper transform of $X$ in the
blowup of $\P$ along $\L$. Since $L \subset X$ and $X$ is
nonsingular this is the blowup of $X$ along $L$.

Thus,
$f_1|_{\G \setminus E}:\G \setminus E \to X \setminus L$
is an isomorphism, $f_1^{-1}(L)= E$ and the
projection $f_1|_E : E \to L$ makes $E$ into a $\bold P^{n-m-1}$ bundle
over $L$. In particular,  for each $p \in L$,
$f_1^{-1}(p)$ is naturally identified with $\bold P(N_p(L/X))$
where $N_p(L/X) = T_pX/ T_pL$ is the normal space of $L$ in $X$ at $p$.
With this identification $f_2|_{f_1^{-1}(p)}$ is just the
linear inclusion $\bold P(N_p(L/X)) \hookrightarrow \bold P(V / V_0)$.

Pick a linear subspace $V_1$ of $V$ complementary to $V$ i.e. so that
$V = V_0 \oplus V_1$ and $\text{dim}(V_1)=N-m$. Then
the projection can be described as follows: if $p \in X \setminus L =
\G \setminus E$ then
$\pi_L(p)$ is the unique point of intersection of $<L,p>$ with $\bold P(V_1)$.
If $p \in E$, then let $q = f_1(p)$ and let $v \in \bold P(N_q(L/X)) \subset
\bold P(V_1)$ be the normal direction at $q$ determined by $p$.
Then $<L,v>$ is a well defined $m+1$ dimensional linear subspace of
$\bold P(V)$ and $f_2(p)$ is the unique point of intersection of
$<L,v>$ with $\bold P(V_1)$.

Now let $X^n \subset \P = \bold P(V)$ be a Roth variety with
$S = S_{0,0,a_1, \dots, a_{n-1}}$ the rational normal scroll having
singular locus $L$ such that $L \subset X \subset S$. Desingularize
the scroll $S$ as in the above construction to obtain the diagram
$$
\CD
\St @>\pi_2>> S \subset \P \\
@V\pi_1VV @. \\
\PP @.
\endCD
$$
Let $\Xt = \overline{\pi_2^{-1}(X \setminus L)}$.

\proclaim{Claim 3.6} If $S \ne S_{0,0,1,1}$ then $\pi_2|_{\Xt}:\Xt \to X$
is an isomorphism. If $S = S_{0,0,1,1}$ there are two possible ways
to desingularize $S$ by the above construction corresponding to the
two rulings of $S_{1,1} \cong \PP \times \PP$. For one of them, $\Xt$ is
isomorphic to $X$ and for the other $\Xt$ is the blow up of
$X$ along $L$.
\endproclaim

First we recall the statement of the version of Zariski's main theorem
that is used in the proof of the claim:
\proclaim{Theorem 3.7} \cite{Mum2, Section III.9, Proposition 1}
Let $X$ be an $n$-dimensional factorial variety and let
$f:X^{\prime} \to X$ be a birational morphism. Then there is a non-empty
open set $U \subset X$ such that
\roster
\item $f|_{f^{-1}(U)}:f^{-1}(U) \to U$ is an isomorphism,
\item If $E_1, \dots, E_k$ are the components of $X^{\prime}
\setminus f^{-1}(U)$ then $\text{dim}(E_i) = n-1$ for all $i$, while
$\text{dim}(\overline{f(E_i)}) \le n-2$.
\endroster
\endproclaim

\demo{Proof of Claim}
Identify $\bold P(V/V_0)$ with $\bold P^{N-2} =
\bold P(V_1 \oplus \dots \oplus V_{n-1}) \subset \bold P(V) = \bold P^N$.
Now let $W \subset \PP \times \bold P^{N-2}$ be defined by:
$$ W = \{(p, <\phi_1(p), \dots , \phi_{n-1}(p)>) \mid p \in \PP \}$$
and consider the diagram of natural projections
$$
\CD
W @>g_2>> S_{a_1,\dots,a_{n-1}} \subset \bold P^{N-2} \\
@Vg_1VV @. \\
\PP @.
\endCD
$$
$g_2$ is an isomorphism and $g_1$ makes $W$ into a $\bold P^{n-2}$-bundle
over $\PP$.

Let $\pi_L:X - - \to \bold P^{N-2}$ be the rational map given by projection
from $L$. Since $X \subset S$, $\pi_L(X) \subset S_{a_1,\dots, a_{n-1}}$ and
thus if $\G$ is the graph of $\pi_L$, then $\G \subset X \times S_{a_1,\dots,
a_{n-1}}$ with diagram of projections
$$
\CD
\G @>f_2>> S_{a_1,\dots,a_{n-1}} \subset \bold P^{N-2} \\
@Vf_1VV @. \\
X @.
\endCD
$$
We obtain a morphism $h:\G \to \PP \times \P$ via the composition
$$
\G \subset X \times S_{a_1,\dots,a_{n-1}} @>(id,g_1\circ g_2^{-1})>>
X \times \PP @>\text{interchange}>> \PP \times X \subset \PP \times \P
$$
such that $h(\G) = \Xt$ and the diagram
$$
\CD
\G @>f_1>> X \\
@VhVV @| \\
\Xt @>\pi_2|_{\Xt}>> X
\endCD
$$
is commutative.

Now $\G\setminus E \cong X \setminus L \cong \Xt \setminus \PP \times L$.
If $p \in L$ then $f_1^{-1}(p) = p \times \bold P(N_p(L/X))$. Then
$h(f_1^{-1}(p))$ is a point $\iff g_1 \circ g_2^{-1}(\bold P(N_p(L/X)))$
is a point $\iff g_2^{-1}(\bold P(N_p(L/X)))$ is a fiber of the
$\bold P^{n-2}$ bundle $g_1:W \to \PP$. But any $n-2$ dimensional
linear space on a nonsingular $n-1$ dimensional scroll must be a fiber as
above unless the scroll is $S_{1,1}$ in which case it could also be a fiber
for the other possible way of making $S_{1,1}$ into a $\PP$ bundle via
the above construction. Thus if $S \ne S_{0,0,1,1}$ then
$\pi_2|_{\Xt}: \Xt \to X$ is a birational bijective morphism and is thus
an isomorphism by Zariski's main theorem. If $S=S_{0,0,1,1}$ then by
Zariski's main theorem either $\pi_2|_{\Xt}^{-1}(L) = \PP \times L$ or
$\pi_2|_{\Xt}:\Xt \to X$ is bijective and thus an isomorphism. In the
first case, $h:\G \to \Xt$ is a bijective morphism. It is easy to see that
$h$ is an immersion and so it is an isomorphism. Thus $\Xt$ is the
blow up of $X$ along $L$. In either case we see that
$\bold P(N_p(L/X)) \cdot \bold P(N_q(L/X)) = 0$ in $A(S_{1,1})$ for $p,q
\in L$ and thus switching to the other desingularization for $S$ we switch
between the two possibilities for $\Xt$.
\qed
\enddemo

By the claim we can assume that $\pi_2|_{\Xt}: \Xt \to X$ is an
isomorphism. In particular, $\Xt \cdot C = 1$. Now $\Xt = \a H + \b F$ for
some integers $\a$ and $\b$. Then,
$$
\align
\Xt \cdot C &= (\a H + \b F) \cdot (H^n - d_SH^{n-1}\cdot F) \\
 &= \a H^{n+1} + (\b - \a d_S) H^{n} \cdot F = \b
\endalign
$$
and so $\b = 1$. Finally, $d = \text{deg}(X) = H^n \cdot \Xt = \a d_S + 1$
and so $\a = \frac{d-1}{N-n}$. To sum up we have proved:

\proclaim{Theorem 3.8}
Let $S^{n+1} \subset \P$ be a rational normal scroll
$S_{0,0,a_1,\dots,a_{n-1}}$ with all $a_i \ge 1$ and $\sum_{i=1}^{n-1}a_i =
N-n$ and let $L$ be the singular locus of $S$. Let
$\Cal E = \Cal O_{\PP} \oplus \Cal O_{\PP} \oplus \Cal
O_{\PP}(a_1) \oplus \dots \oplus \Cal O_{\PP}(a_{n-1})$ and let
$\pi_1:\bold P(\Cal E^*) \to \PP$ be a $\bold P^n$ bundle with diagram:
$$
\CD
\bold P(\Cal E^*) @>\pi_2>> S \subset \P \\
@V\pi_1VV @. \\
\PP @.
\endCD
$$
that desingularizes $S$ such that the morphism $\pi_2$ is given by the
complete linear series $|\Cal O_{\bold P(\Cal E^*)}(1)|$. Let H be the
pullback of a hyperplane in $\P$ via $\pi_2$ and let $F$ be a fiber of $\pi_1$.
Then for every $b > 0$ there exists a nonsingular variety $\Xt \in |bH +F|$
such that if $X = \pi_2(\Xt)$ then $\pi_2|_{\Xt}:\Xt \to X$ is an
isomorphism and $X$ is a linearly normal Roth variety of
degree $b(N-n)+1$ (with associated scroll $S$).
Conversely, if $X$ is a Roth variety of degree $d$ such that
$L \subset X \subset S$, then $b = \frac{d-1}{N-n}$ is an integer and there
is a desingularization of $S$ as above; $\pi_2: \bold P(\Cal E^*) \to S$ and
$\Xt \in |bH + F|$ which is mapped isomorphically onto $X$ by $\pi_2$.
(In particular, Roth varieties are linearly normal).
\endproclaim
Note that further properties of Roth varieties are given in Proposition 3.18.

Now let $d = \text{deg}(X)$. Then $d = H^n \cdot \Xt = H^n \cdot (bH+F)
= bd_S+1$. Let $\CX$ denote a double-point divisor of $X \subset \P$.

\proclaim{Proposition 3.9}
$\CX \cdot L = 0$. Thus $|\CX|$ is not ample.
\endproclaim

\demo{Proof}
$\CX \cdot L = ((d-n-2)H - K_X) \cdot L$.
Since $\Xt$ is isomorphic to $X$ we can write the intersection
occuring in $A(X)$ as an intersection occuring in $A(\Xt)$. Thus,
$$\CX \cdot L = ((d-n-2)H- K_{\Xt}) \cdot (\pi_2|_{\Xt})^{-1}(L).$$
$(\pi_2|_{\Xt})^{-1}(L)$
is rationally equivalent in $A(\St)$ to $(\PP \times L) \cdot \Xt$.
We can then use the
projection formula and adjunction to write the intersection in $A(\Xt)$
as an intersection in $A(\St)$. Thus,
$$\CX \cdot L = ((d-n-2)H - K_{\St} - \Xt) \cdot (\PP \times L) \cdot \Xt.$$
Now to verify the theorem it suffices
to directly evaluate the right hand side in the above equation using the
intersection theory for $\St$ developed earlier. Explicitly,
$$
\align
\CX \cdot L &= ((d-n-2)H-K_{\St}-\Xt) \cdot (\PP \times L) \cdot \Xt \\
&=((d-n-2)H-(-(n+1)H + (d_S-2)F) - (b H+F)) \\
& \qquad  \cdot (H^{n-1} - d_S H^{n-2} \cdot F) \cdot (b H+F) \\
&=(b(d_S-1)H + (1-d_S)F)\cdot (bH^n + (1-b d_S)H^{n-1}\cdot F) \\
&=b^2(d_S-1)H^{n+1} + b(d_S-1)(1-b d_S)H^n \cdot F + b(1-d_S)H^n \cdot F \\
&=0.
\qed
\endalign
$$
\enddemo

Let $\pi$ be the {\sl sectional genus} of $X$.
A generic $(N-n+1)$-plane section
of $X$ is a smooth curve and $\pi$ is by definition the genus of this
curve. $\pi$ can be computed by the formula:
$$2\pi - 2 = K_X \cdot H^{n-1} + (n-1)H^n.$$

\proclaim{Proposition 3.10}
$$\pi = \frac{1}{2} \frac{(d-1)(d-(N-n+1))}{(N-n)}$$
\endproclaim

\demo{Proof}
First recall that $d = b d_S +1$ and $d_S = N-n$ so that $d = b(N-n)+1$.
Now use the formula $2\pi - 2 = K_X \cdot H^{n-1} + (n-1)H^n.$ As in
the previous proof, first express the right hand side as an intersection
in $A(\Xt)$ and then use adjunction and the projection formula to write it
as an intersection in $A(\St)$. In detail,
$$
\align
2\pi - 2 &= K_{\Xt} \cdot H^{n-1} + (n-1)H^n \\
&=(K_{\St}+\Xt) \cdot \Xt \cdot H^{n-1} + (n-1)H^n \cdot \Xt \\
&=((-(n+1)H+(d_S-2)F)+(bH+F))\cdot (bH+F)\cdot H^{n-1} \\
& \qquad + (n-1)H^n \cdot(bH+F) \\
&=((b-(n+1))H + (d_S-1)F)\cdot (bH^n + H^{n-1} \cdot F) \\
& \qquad + (n-1)(bH^{n+1}+H^n \cdot F) \\
&=((b-(n+1))b + (n-1)b)H^{n+1}  \\
& \qquad + ((b-(n+1))+(d_S-1)b + (n-1))H^n \cdot F \\
&=b^2d_S-bd_S -2
\endalign
$$
Thus,
$$
\pi = \frac{1}{2}bd_S(d_S-1) = \frac{1}{2} \frac{(d-1)(d-(N-n+1))}{(N-n)}
\qed
$$
\enddemo

At this stage it is convenient to state and prove some facts about
normality.

\proclaim{Proposition 3.11}
Let $X^{n} \subset \P$ be a  nondegenerate  variety
and let $Y^{n-1} \subset \bold P^{N-1}$ be an irreducible
hyperplane section of $X$.
\roster
\item If $X$ is $k-1$ normal and $Y$ is $k$ normal then $X$ is $k$ normal.
\item If $Y$ is linearly normal then $X$ is linearly normal.
\item If $Y$ is projectively normal then $X$ is projectively normal.
\item If $X$ is $k$ normal then:
$$
Y \text{is\ } k \text{\ normal} \iff H^1(X, \Cal O_X(k-1)) \hookrightarrow
H^1(X, \Cal O_X(k)).
$$
\item If $Y$ is projectively normal then $H^1(X, \Cal O_X(k)) = 0$ for all
$k \ge 0$.
\endroster
\endproclaim

\demo{Proof}
For \tri1 we have a commutative diagram with exact rows:
$$
\matrix
H^0(\P, \Cal O_{\P}(k-1)) & \hookrightarrow & H^0(\P, \Cal O_{\P}(k)) & \to &
H^0(\bold P^{N-1}, \Cal O_{\bold P^{N-1}}(k)) & \to 0 \\
\downarrow \alpha & & \downarrow \beta & & \downarrow \gamma & \\
H^0(X, \Cal O_X(k-1)) & \hookrightarrow & H^0(X, \Cal O_X(k)) & @>\phi>> &
H^0(Y, \Cal O_Y(k)) &
\endmatrix
$$
Hence if $X$ is $k-1$ normal and $Y$ is $k$ normal then $\alpha$ and
$\gamma$ are surjective and thus by the Snake lemma, $\beta$ is surjective
i.e. $X$ is $k$ normal.

For \tri2, apply \tri1 observing that $X$ is nondegenerate means that
$X$ is $0$ normal. For \tri3, apply \tri1 inductively.
For \tri4, if $\beta$ is surjective then $\text{image}(\phi) = \text{image}
(\gamma)$. Thus $Y$ is $k$ normal $\iff \gamma$ is surjective $\iff \phi$
is surjective $\iff H^1(X, \Cal O_X(k-1)) \hookrightarrow H^1(X, \Cal O_X(k))$.
By \tri4 and \tri3 we have
$$
H^1(X, \Cal O_X) \hookrightarrow H^1(X, \Cal O_X(1)) \hookrightarrow
H^1(X, \Cal O_X(2)) \hookrightarrow \dots \hookrightarrow H^1(X,\Cal O_X(k))
$$
Thus by Serre vanishing \tri5 follows.
\qed
\enddemo

\proclaim{Proposition 3.12}
Let $X^n \subset \P$ be a nonsingular, nondegenerate projective variety
and suppose that the generic $N-n+1$ plane section of $X$ is a nonsingular
projectively normal curve. Then $H^i(X, \Cal O_X(k))= 0$ for  $1 \le i \le
n-1$ and all $k \ge 0$.
\endproclaim

\demo{Proof}
By induction on $n$. For $n=1$ there is nothing to prove. So suppose the
result is true for all varieties of dimension $< n$. Let $Y$ be a
generic hyperplane section of $X$. Then
$H^i(Y, \Cal O_Y(k)) = 0$ for $1 \le i \le n-2$ by the induction hypothesis.
The long exact sequence of cohomology associated to the exact sequence
$$
0 \to \Cal O_X(k-1) \to \Cal O_X(k) \to \Cal O_Y(k) \to 0
$$
shows that $H^i(X, \Cal O_X(k-1)) \hookrightarrow H^i(X, \Cal O_X(k))$
for $2 \le i \le n-1$ and $k \ge 0$.
 Thus, by Serre vanishing $H^i(X, \Cal O_X(k)) = 0$
for $2 \le i \le n-1$. By Proposition 3.11 \tri3 and \tri5,
$H^1(X, \Cal O_X(k)) = 0$ for all $k \ge 0$. Finally, vanishing for all
$k \in \bold Z$ follows by Kodaira's vanishing theorem.
\qed
\enddemo

\subheading{Remark 3.13} Let $X^n$ be a  nonsingular
variety such that $h^1(X, \Cal O_X) \ge 1$ and let
$L$ be a very ample line bundle on X.
Let $Y_k \in |L^{\otimes k}|$ be a generic element. We have the exact
sequence:
$$ 0 \to \Cal O_X \to L^{\otimes k} \to L^{\otimes k}|_{Y_k}
\to 0 $$
Which by Serre vanishing gives that
$$ h^0(Y_k,
L^{\otimes k}|_{Y_k}) - h^0(X, L^{\otimes k}) = h^1(X, \Cal O_X)-1
\ge 0
$$
for $k \gg 0$. Now assume $k \gg 0$. Let $\phi: X \to \bold P(H^0(X,
L^{\otimes k})) = \P$ be the morphism determined by $|L^{\otimes k}|$. Let
$\phi(X) = X^{\prime}$ and $\phi(Y_k) = Y_k^{\prime}$. Then $X^{\prime}$ is
linearly normal, $Y_k^{\prime}$ is a nonsingular hyperplane section of
$X^{\prime}$ but by the above inequality, $Y_k^{\prime}$ is not linearly
normal. Thus the converse of \tri2 above is false.

The above construction produces examples of high codimension. One can ask
for conditions on $X^n \subset \P$ of small codimension so that a hyperplane
section of $X$ is linearly normal. For example, C. Peskine's problem: find
$d$ such that a smooth surface in $\bold P^4$ of degree greater than $d$ has
linearly normal hyperplane sections.

\subheading{Remark 3.14} The converse to \tri3 is also false. See
Example 5.8 for a counterexample.

Let $X^n \subset \P$ be a nondegenerate variety of degree $d$ and
$\phi: X^{\prime} \to X$ a resolution of $X$ i.e. $X^{\prime}$ is
nonsingular and $\phi$ is a birational morphism. Such a resolution
exists by Hironaka's resolution of singularities. Then the geometric
genus of $X$, $p_g(X) = h^n(X^{\prime}, \Cal O_{X^{\prime}})$ is a
well defined invariant of $X$.  Let $M = \fwd[]{d-1}{N-n}$ and let
$\epsilon$ be defined by $d-1 = M(N-n) + \epsilon$. Then there is a
bound due to J. Harris on the geometric genus of $X$ in terms of $d$,
$n$, and $N$:
\proclaim{Theorem 3.15} \cite{Ha1, page 44}
$$
p_g(X) \le \binom{M}{n+1}(N-n) + \binom{M}{n}\epsilon
$$
where the binomial coefficient $\binom ab = 0$ when $b > a$.
\endproclaim

This theorem motivates Harris' definition of Castelnuovo variety:
\subheading{Definition 3.16}
$X$ is a {\sl Castelnuovo variety} if
\roster
\item $$ p_g(X) = \binom{M}{n+1}(N-n) + \binom{M}{n}\epsilon$$
\item $d \ge n(N-n) + 2$.
\endroster

The second condition in the definition is intended to ensure that
$p_g(X) \ne 0$.
If $n=1$, Harris' bound is just the classical Castelnuovo bound for the
genus of a space curve. The curves whose genus achieve this bound are
Castelnuovo curves if $d \ge (N-n) + 2$ and rational normal curves
if $d = (N-n) + 1$. In either case, these curves are projectively
normal.

Now let $X^n \subset \P$ be a Roth variety of degree $d = b(N-n)+1$
contained in a rational normal scroll $S = S_{0,0,a_1,\dots,a_{n-1}}$ with
all $a_i \ge 1$ and singular locus $L$ such that $L \subset X \subset S$.
Then by Proposition 3.10, the geometric genus of $X$ is
$$ \pi = \frac{1}{2} \frac{(d-1)(d-(N-n+1))}{(N-n)}$$
so a generic $\N' = N-n+1$ plane section of X is a
 curve $Y$ of degree $d$ and
genus
$$g(Y) = \frac{1}{2} \frac{(d-1)(d-\N')}{(\N' -1)}$$
in $\bold P^{\N'}$. In this case, the Castelnuovo bound for the genus of
$Y$ is
$$g(Y) \le \binom{M}{2} (\N' - 1) + M \epsilon.$$
Now $M = \frac{d-1}{\N' - 1} = b$ and so $\epsilon = 0$. Thus,
$$
g(Y) \le \frac{1}{2} \frac{d-1}{\N' -1} \left({\frac{d-1}{\N' -1} - 1}\right)
(\N' - 1) =\frac{1}{2} \frac{(d-1)(d-\N')}{(\N' - 1)}
$$

If $b \ge 2$ then we've shown that a generic $N-n+1$ plane section
of $X$ is a Castelnuovo curve. If $b=1$ then $d = (N-n)+1$ and so by
classification of varieties of minimal degree, $X$ is a rational
normal scroll. (In fact, by Proposition 5.9,
$N_{L/X} = \Cal O_{L}(1-a_1) \oplus \dots \oplus \Cal O_L(1-a_{n-1})$.
It follows by Proposition 5.10 that
$X$ is a rational normal scroll $S_{1,a_1, \dots , a_{n-1}}$).
In either case, by Proposition 3.11 \tri3, $X$ is projectively normal.
By Proposition 3.12, $H^i(X, \Cal O_X(k)) = 0$ for $1 \le i \le n-1$ and
$k \ge 0$.

Harris states that a nonsingular nondegenerate variety $X^n \subset \P$ is
Castelnuovo if and only if its generic hyperplane section is.
The ``if'' direction of this theorem is false, essentially because the
condition \tri2 in the definition of Castelnuovo variety is not
stable under taking hyperplane sections.
For example, if $X \subset \P$ is a Roth surface of degree $d = 2(N-n)+1$
then $X$ is not a Castelnuovo variety but a generic hyperplane section
of $X$ is a Castelnuovo curve.
However, provided that the degree of $X$ is in the required range then
Harris' result is correct. Namely,
\proclaim{Proposition 3.17}
If $X^n \subset \P$ is a nonsingular nondegenerate variety of
 degree $d \ge n(N-n)+2$
and $Y$ is a generic hyperplane section of $X$ then $Y$ is Castelnuovo
if and only if $X$ is.
\endproclaim
Thus, if $X$ is a Roth variety of degree $d = b(N-n)+1$ with $b \ge n+1$
then $X$ is a Castelnuovo variety.

Probably the two most characteristic properties of Castelnuovo varieties
are:
\roster
\item a generic $N-n+1$ section of $X$ is a Castelnuovo curve
\item $X$ is a divisor in a rational normal scroll cut out by the
      quadrics containing $X$.
\endroster
In fact, Harris proves \tri2 for Castelnuovo varieties by essentially
showing that \tri1 implies \tri2. Hence one may be tempted to
take \tri1 as the definition of a Castelnuovo variety. With this
definition, a Roth variety that is not a rational normal scroll
is a Castelnuovo variety.  Later we will give  another very geometric proof
of \tri2 for Roth varieties that are not rational normal scrolls.

To sum up, we've proved:
\proclaim{Theorem 3.18}
Let $X^n \subset \P$ be a Roth variety of degree $d = b(N-n)+1$
and $S$ a rational normal scroll $S_{0,0,a_1,\dots, a_{n-1}}$ with
all $a_i \ge 1$, $\sum_{i=1}^{n-1}a_i = N-n$, and singular locus
$L$ such that $L \subset X \subset S$.
\roster
\item If $b = 1$ then $X$ is a rational normal scroll
      $S_{1,a_1,\dots, a_{n-1}}$.
\item If $b \ge 2$ then the generic $N-n+1$ plane section of $X$ is
	a Castelnuovo curve.
\item $X$ is projectively normal.
\item $H^i(X, \Cal O_X(k)) = 0$ for $1\le i \le n-1$ and all $k \in \bold Z$.
\item If $b \ge n + 1$ then $X$ is a Castelnuovo variety.
\endroster
\endproclaim
Note that \tri3 and \tri4 imply that $X$ is arithmetically
Cohen-Macaulay.

\newpage

\heading 4. Proof of the main reduction \endheading

\proclaim{Proposition 4.1} Let $X^n \subset \bold P^N$ be a nonsingular,
nondegenerate, linearly normal variety of dimension $n \ge 2$ and
codimension $N-n \ge 2$. Let $p \ne q \in X$ be such that
$S(p,X) = S(q,X)$ and let $L = <p,q>$. Then,
\roster
\item $L \subset X$.
\item $S(L,X)$ is a rational normal scroll $S_{0,0,a_1,\dots,a_{n-1}}$
with all $a_i \ge 1$ and \break $\sum_{i=1}^{n-1}a_i = N-n$.
\item $L$ is the singular locus of $S(L,X)$.
\item $S(r,X) = S(L,X)$ for all $r \in L$.
\endroster
In particular, $X$ is a Roth variety.
\endproclaim

The main corollary is a classification of those varieties whose
double-point divisor determines an ample linear system:
\proclaim{Theorem 4.2} Let $X^n \subset \P$ be a nonsingular, nondegenerate
 variety of codimension $N-n \ge 2$. Let $|\CX|$ denote the
linear system determined by the double-point divisor on $X$.
Then the following are equivalent:
\roster
\item $|\CX|$ is ample.
\item $|\CX|$ separates points.
\item $X$ is not an isomorphic projection of a Roth variety.
\endroster
If, in addition, $X$ is linearly normal then \tri3 can be replaced by:
\roster
\item"$(3^{\prime})$" $X$ is not a Roth variety.
\endroster
\endproclaim

\demo{Proof}
By Remark 2.1 it suffices to prove the theorem when $X$ is linearly normal.
Then \tri3 $\implies$ \tri2 follows by Corollary 2.5 and Proposition 4.1.
\tri2 $\implies$ \tri1 follows by Corollary 2.7 and \tri1 $\implies$
\tri3 follows by Proposition 3.9.
\qed
\enddemo

An alternative version of this theorem can be given which amounts to
a classification of those varieties with certain highly degenerate
associated secant varieties.
\proclaim{Theorem 4.3} Let $X^n \subset \P$ be a nonsingular, nondegenerate,
variety of codimension $N-n \ge 2$.
Then the following are equivalent:
\roster
\item There exists a positive dimensional subvariety
$W$ of $X$ such that the secant variety $S(W,X)$ has dimension $n+1$.
\item There is $p \ne q \in X$ such that $S(p,X) = S(q,X)$.
\item $X$ is an isomorphic projection of a Roth variety.
\endroster
If, in addition, $X$ is linearly normal then \tri3 can be replaced by:
\roster
\item"$(3^{\prime})$" $X$ is a Roth variety.
\endroster
\endproclaim

\demo{Proof}
Clearly it suffices to prove the theorem in the case when $X$ is linearly
normal. \tri1 $\implies$ \tri2 is clear. \tri2 $\implies$ \tri3 is
Proposition 4.1. \tri3 $\implies$ \tri1 since if $X$ is a Roth variety
with $L \subset X \subset S$ then for all $p \in L$, $S(p,X) = S$ so
that $S(L,X)=S$ is $n+1$ dimensional.
\qed
\enddemo

To begin the proof of Theorem 4.1, we quote a Bertini-like proposition that
will be used often:
\proclaim{Proposition 4.4} \cite{Ha2, Proposition 18.10}
Let $X^n \subset \P$ be a nondegenerate variety
of dimension $n \ge 1$ and let $H$ by a hyperplane intersecting $X$
generically transversally i.e. for a generic point $p \in X \cap H,\  \Tt_pX
\not\subset H$. Then if $Y = X \cap H$ then $Y$ is nondegenerate as a
(possibly reducible) subvariety of $H \cong \bold P^{N-1}$.
\endproclaim

We also need a corollary of the Fulton-Hansen connectedness theorem:
\proclaim{Theorem 4.5} \cite{FL, Theorem 4.1}
Let $X^n \subset \P$ be a variety of dimension $n \ge 2$ and let $H$ be
a hyperplane in $\P$. Then $H \cap X$ is connected.
\endproclaim

\proclaim{Proposition 4.6} Let $X^n \subset \P$ be a nonsingular,
nondegenerate variety of dimension $n \ge 2$ and codimension $N-n \ge 2$.
Let $Z^{n-1} \subset \bold P^{N-2} \subset \P$ be a variety and
let $L$ be a line in $\P$ disjoint from $\bold P^{N-2}$. Let $\CLZ$ be the
cone over $Z$ with vertex $L$. $\CLZ$ is an $n+1$ dimensional variety. Suppose
that $X \subset \CLZ$. Then $L \subset X$.
\endproclaim
\demo{Proof}
The proof is by induction on $n$. Suppose that $n = 2$ and that
$L \not\subset X$ so that $L \cap X = \{p_1, \dots, p_j \}$ is a
(possibly empty) finite set of points. Let $q_i \in \Tt_{p_i}X \setminus L$.
Let $H$ be a generic hyperplane containing $L$ such that $q_i \notin H$ for
all $i$. Let $Y = X \cap H$. By Bertini's theorem, $Y$ is nonsingular
except possibly at $p_1, \dots, p_k$. Since $\Tt_{p_i}X \not\subset H$ it
follows by the Jacobian criterion that $Y$ is nonsingular at $p_i$ for all
$i$. Hence $Y$ is nonsingular. By Proposition $4.4$, $Y \subset H$ is
nondegenerate.
By the Fulton-Hansen connectedness theorem, $Y$ is connected. Hence $Y$ is
irreducible. On the other hand, if $H \cap Z = \{q_1, \dots, q_k \}$ then
since $X \subset \CLZ$, $Y \subset H \cap \CLZ =
<L,q_1> \cup \dots \cup <L,q_k>$.
Since $Y$ is irreducible, $Y \subset <L,q_i>$ for some $i$. This
contradicts the nondegeneracy of $Y \subset H$.

If $n \ge 3$ and $L \not\subset X$ then by slicing with a generic hyperplane
$H$ through $L$ we obtain $Y = H \cap X$. By the same argument as in the
$n=2$ case  it follows that $Y$ is nonsingular and a nondegenerate subvariety
of $H$. By the induction hypothesis $L \subset Y$, which is a contradiction.
Hence $L \subset X$.
\qed
\enddemo

\proclaim{Proposition 4.7} Let $X^n \subset \P$ be a nonsingular,
nondegenerate variety of dimension $n \ge 2$ and codimension $N-n \ge 2$.
Let $p \ne q \in X$ be such that $S(p,X) = S(q,X)$. Let $L = <p,q>$. Then
$L \subset X$ and $S(r,X) = S(L,X)$ for all $r \in L$.
\endproclaim

\demo{Proof}
Let $S = S(p,X)$ and let $\bold P^{N-2}$ be a generic $N-2$ plane in $\P$.
Then $\bold P^{N-2} \cap L = \emptyset$ and $Z := \bold P^{N-2} \cap S$
is an $n-1$ dimensional variety. Applying Proposition 2.10 to
$W = \{p,q \}$ we see that $W \subset \Tt_pX \cap \Tt_qX$. Thus
$L \subset \Tt_pX \subset S$. Let $r \in Z$. We claim that $<L,r> \subset S$.
Since $r \in S$, either
$r \in \Tt_pX$ or the line $<p,r>$ is a proper bisecant of $X$. In the
first case, $<L,r> \subset \Tt_pX \subset S$. If $r \in \Tt_qX$ then also
$<L,r> \subset S$. So suppose that $r \notin \Tt_qX$ and that $<p,r>$ is
a proper bisecant of $X$. Then since $S(p,X) = S(q,X)$, for generic $w \in
<p,r>$, the line $<w,q>$ is a proper bisecant line of $X$. Varying $w$ we
obtain a curve $C \subset <L,r>$ such that $C \ne L$. It follows that
$<L,r> \subset S$.
Let $T$ be the cone over $Z$ with vertex $L$. Then $T \subset S$ and
since $\text{dim}(Z) = \text{dim}(S)$ and $S$ is irreducible it
follows that $T = S$. Thus, by Proposition 4.6, $L \subset X$. Clearly
$S(r,X) = S(L,X)$ for all $r \in L$.
\qed
\enddemo

\subheading{Set up of notation 4.8}
Now let $X^n \subset \P$ be a nonsingular, nondegenerate variety of
dimension $n \ge 2$ and codimension $N-n \ge 2$. Let $L^m \subset X$ be
a possibly reducible subvariety of $X$ such that $m \ge 1$
and $S(p,X)= S(L,X)$ for all $p \in L$. By Proposition 4.7 we can assume
that $L$ is a linear space. Let $\bold P(V_1)$ be a generic
$N-m-1$ plane in $\bold P^N$. Then $\bold P(V_1)
\cap L = \emptyset$ and $Z := \bold P(V_1) \cap S(L,X)$ is an $n-m$
dimensional variety.  $Z \subset \bold P(V_1)$ is
nondegenerate by Proposition 4.4.  Let $\CLZ$ be the cone over $Z$ with
vertex $L$.
\proclaim{Proposition 4.9} With notation as in 4.8,
\roster
\item If $q \in S(L,X)$ then $<L,q> \subset S(L,X)$.
\item $\CLZ = S(L,X)$.
\endroster
\endproclaim

\demo{Proof}
Let $q \in S(L,X)$. Then by Proposition 2.9,
$S(L,X) = \bigcup_{p \in L} \Tt_pX$, so there is $p \in L$ such
that $q \in \Tt_pX$. It follows that
$<L,q> \subset \Tt_pX \subset S(L,X)$. Since $Z \subset S(L,X)$,
it follows that $\CLZ \subset S(L,X)$. On the other hand, $\CLZ$
and $S(L,X)$ are both irreducible varieties of dimension
$n+1$. Hence $\CLZ = S(L,X)$.
\qed
\enddemo

\proclaim{Proposition 4.10}
With notation as in 4.8, suppose that $n-m \ge 2$.
Let $H$ be a generic hyperplane containing $L$. Then $Y = X \cap H$ is
a nonsingular, nondegenerate variety and $S(p,Y)= S(L,Y)$ for
all $p \in L$.
\endproclaim

\demo{Proof}
Since $L \subset \bigcap_{p\in L}\Tt_pX$ by Proposition 2.10  and since $L$
and $\bold P(V_1)$
are complementary linear spaces in $\P$,
$\Tt_pX \cap \bold P(V_1)$ is an
$n-m-1$ plane for every $p \in L$.
Moreover, since $S(L,X) = \bigcup_{p \in L} \Tt_pX$,
this plane is contained in $Z$.
Let $H$ be a generic hyperplane containing $L$. By Bertini's theorem,
$Y = X \cap H$ is a possibly reducible variety that is nonsingular off of $L$.
By Proposition 4.4, $Y$ is nondegenerate. Moreover, by the Jacobian criterion,
if $p \in L$ then $Y$ is nonsingular at $p$ if
$\Tt_pX  \not\subset H$ which holds if
$\Tt_pX \cap \bold P(V_1) \not\subset H \cap \bold P(V_1)$.
But $H \cap \bold P(V_1)$ is a generic hyperplane in $\bold P(V_1)$ which
meets $Z$ in an $n-m-1 \ge 1$ dimensional variety
 that is not a linear space since
$Z$ is not a linear space. Hence the $n-m-1$ plane $\Tt_pX \cap \bold P(V_1)$
is not contained in $H \cap \bold P(V_1)$ for any $p \in L$. Thus
$Y$ is nonsingular. Now by the Fulton-Hansen connectedness theorem,
$Y$ is connected and thus $Y$ is irreducible.
$S(p,Y)=S(L,Y)$ for all $p \in L$ since $S(p,X) \cap H = S(p, Y)$.
\qed
\enddemo

Let $X^n \subset \P$ be a nonsingular, nondegenerate projective
variety. Let $L^m$ be a linear subspace contained in $X$ with $m \ge 1$ and
let $\bold P(V_1)$ be an $n-m-1$ plane in $\P$ that does not meet $L$. Let
$\pi_L:X - - \to  \bold P(V_1)$
be the rational map determined by projection
from the center $L$ onto $\bold P(V_1)$. Let $\G \subset X \times \bold P(V_1)$
be the graph of $\pi_L$. $\G$ is also  the blow up of $X$ along $L$. Let $E$
be the exceptional divisor of the blow up. We have the diagram:
$$
\matrix
E & \subset & \G & @>f_2>> & \bold P(V_1) \\
\Big\downarrow &  & f_1\Big\downarrow \  &  & \Big\Vert \\
L & \subset & X & -- \to & \bold P(V_1)
\endmatrix
$$

\proclaim{Proposition 4.11}
With notation as in 4.8:
\roster
\item $Z$ is an $n-m$ dimensional variety, not a linear space, and
nondegenerate as a subvariety of $\bold P(V_1)$.
\item $Z = \pi_L(X)$.
\item $f_2|_{E} : E \to Z$ is surjective.
\endroster
\endproclaim

\demo{Proof}
\tri1 is part of the set up of 4.8. By Proposition 4.9 \tri1,
 $\pi_L(X) \subset Z$. Now, since
$S(L,X) = \cup_{p \in L} \Tt_pX$,
$$Z = S(L,X) \cap \bold P(V_1) = \cup_{p \in L}(\Tt_pX \cap \bold P(V_1))
= \cup_{p \in L}f_2(f_1^{-1}(p)) = f_2(E).$$
Thus, \tri2 and \tri3 follow.
\qed
\enddemo

We record an elementary and well known lemma used in proving the next
proposition.
\proclaim{Lemma 4.12}
Let $C$ be a curve and $k \ge 2$ an integer. Then there is no nonconstant
morphism $\bold P^k \to C$.
\endproclaim
\demo{Proof}
Suppose that $\phi:\bold P^k \to C$ is a nonconstant morphism. Since $C$ is
a curve, the morphism is surjective. Let $p,q \in C, \ p \ne q$. Then
$\phi^{-1}(p)$ and $\phi^{-1}(q)$ are $k-1$ dimensional subvarieties of
$\bold P^k$ which do not intersect. Since $k \ge 2$ this contradicts
B\'ezout's theorem.
\qed
\enddemo

\proclaim{Proposition 4.13}
With notation as in 4.8, $L$ is a line.
\endproclaim
\demo{Proof}
Using Proposition 4.10, by repeatedly slicing by generic hyperplanes through
$L$
 we can  assume that $n-m=1$. The blow up of $X$ along $L$ is then
isomorphic to $X$ since $L$ is a divisor in $X$. Thus $L = E$ and we
obtain a surjective morphism $L \cong E @>f_2|_{E}>> Z$ from a projective
space of dimension $m$ to the curve $Z$. Hence $m=1$ by the previous lemma.
\qed
\enddemo

As mentioned in the introduction, the next proposition follows the main
line of argument of the proof of \cite{So, Lemma 1.1.1}.
\proclaim{Proposition 4.14} Let $X \subset \bold
P^N$ be a nonsingular, nondegenerate
and linearly normal projective surface with $N \ge 3$. Let $L \subset X$
be a line such that $S(p,X) = S(L,X)$ for all $p \in L$. Then
$S(L,X)$ is a rational normal scroll $S_{0,0,N-2}$.
\endproclaim
\demo{Proof}
The blow up of $X$ along $L$ is isomorphic to $X$ since $L$ is a divisor in
$X$. Thus identifying $\G$ with $X$ and $E$ with $L$ we obtain the diagram:
$$
\matrix
X & @>f_2>> & Z & \subset \bold P^{N-2} \\
\bigcup & & \Big\Vert \\
L & \twoheadrightarrow & Z
\endmatrix
$$
where $Z$ is a curve. Since $Z$ is nondegenerate, the morphism $f_2$ is
the morphism determined by a base point free subsystem $\Lambda$ of the
linear system $|H-L|$ and $\text{dim}(\Lambda)=N-2$. Since X is
nondegenerate and linearly normal, $\text{dim}|H| = N$. Since vanishing
on the line $L$ imposes two conditions on linear forms, $\text{dim}|H-L|
=N-2$. Thus $\Lambda$ is the complete linear system $|H-L|$. Hence $Z$ is
linearly normal. Since $\PP = L \twoheadrightarrow Z$, we have that
$Z \cong \PP$ and that $Z$ is a nondegenerate rational normal curve of
degree $N-2$. Since $S(L,X) = \CLZ$, it follows that $S(L,X)$ is a
rational normal scroll $S_{0,0,N-2}$.
\qed
\enddemo

To prove the next proposition we need to use the classification
of varieties of minimal degree:
\proclaim{Theorem 4.15} \cite{Ha1, page 51}
Let $X^n \subset \P$ be a nondegenerate variety of degree $N-n+1$. Then
$X$ is one of:
\roster
\item a quadric hypersurface.
\item a cone over the Veronese surface in $\bold P^5$.
\item a rational normal scroll.
\endroster
Moreover, $X$ contains an $n-1$ plane $\iff$ $X$ is a rational normal scroll.
\endproclaim

\proclaim{Proposition 4.16}
Let $X^n \subset P^N$ be a nonsingular, nondegenerate and linearly normal
projective variety of codimension  $N -n \ge 2$.
Let $L \subset X$ be a line such that
$S(p,X) = S(L,X)$ for all $p \in L$. Then $S(L,X)$ is a rational normal
scroll $S_{0,0,a_1,\dots, a_{n-1}}$ with all $a_i \ge 1$ and
$N-n = \sum_{i=1}^{n-1}a_i$. Moreover, $L$ is the singular locus of the scroll.
\endproclaim

\demo{Proof}
The proof is by induction on $n$. The case of $n=2$ was done in the last
proposition. So suppose that $n \ge 3$ and that the result is true for
varieties of dimension less that $n$. Take a generic hyperplane $H$
containing $L$ and consider $Y = X \cap H$. By Proposition 4.10
we know that $Y$ is a nonsingular,
 nondegenerate variety with $S(p,Y) = S(L,Y) =
S(L,X) \cap H$ for all $p \in L$. By the induction
hypothesis, $Y$ is an isomorphic projection of a Roth variety. Thus
by Theorem 3.18 \tri4, $H^1(Y, \Cal O_Y) = 0$.
The long exact cohomology sequence associated
to the short exact sequence
$$
0 \to \Cal O_X(-H) \to \Cal O_X \to \Cal O_Y \to 0
$$
gives that $H^1(X, \Cal O_X(-H)) \to H^1(X, \Cal O_X)$ is surjective.
But by Serre duality $H^1(X, \Cal O_X(-H)) \cong H^{n-1}(X, \Cal O_X(
K_X + H))$ which vanishes by Kodaira's vanishing theorem. Hence
$H^1(X, \Cal O_X) = 0$. Thus, by Proposition 3.11 \tri2,
Y is linearly normal and
so $Y$ is a Roth variety and $S(L,Y)$ is a rational normal scroll.
Since $S(L,Y) = S(L,X) \cap H$ the classification of varieties of
minimal degree shows that $S(L,X)$ is a rational normal scroll.

Let $M$ be the singular locus of $S(L,X)$. It remains to verify that
$M = L$. By Proposition 4.9 \tri2 it follows that $L \subset M$.
Since $S(L,X)$ is a rational normal scroll, $S(L,X)$
is a cone with vertex $M$. Pick a line $\Lp$ in $M$. Then
$S(L,X)$ is also a cone over a rational normal scroll $Z^{n-1}$ with vertex
$\Lp$ such that $<Z>$ and $\Lp$ are complementary linear subspaces of
$\P$. Hence by Proposition 4.6, $\Lp \subset X$. Hence $M \subset X$.
Clearly $S(p,X) = S(M,X)$ for all $p \in M$. Hence $M = L$ by
Proposition 4.13.
\qed
\enddemo

\newpage

\heading 5. Applications of the classification and further properties
of Roth varieties \endheading

\proclaim{Theorem 5.1} Let $X^n \subset \P$ be a Roth variety
of degree $d=b(N-n)+1$ and $S$ a rational normal scroll
$S_{0,0.a_1, \dots, a_{n-1}}$ with all $a_i \ge 1$, $\sum_{i=1}^{n-1}a_i
=N-n$, and singular locus $L$ such that $L \subset X \subset S$. Then,
\roster
\item if $b \ge 2$ then the intersection of all quadrics containing $X$
is $S$. Hence the associated scroll $S$ is uniquely determined by $X$.
\item the generic $N-n+1$ plane section of $X$ containing $L$ is
$$L \cup C_1 \cup \dots \cup C_{N-n}$$
where $C_i$ are nonsingular plane curves of degree $b$, $C_i \ne L$ and
if $p_i \in C_i \setminus L$ then $C_i \subset <L,p_i>$. Moreover,
$C_i \cap C_j = \emptyset$ if $i \ne j$ and $C_i \sim C_j$ in $A(X)$.
\item If $n \ge 3$ the generic hyperplane section of $X$ containing
$L$ is a Roth variety.
\endroster
\endproclaim

\demo{Proof}
We will prove \tri2 first. Using Proposition 4.10,
 by taking successive
hyperplane sections through $L$ we can assume that $n=2$. Then
projection from $L$ extends to a morphism $\pi_L: X \to S_{N-2}
\subset \bold P^{N-2}$ as in
the proof of Proposition 4.14. Let $\H'$ be a generic hyperplane
in $\bold P^{N-2}$ and let $H = <L, \H'>$. Suppose that
$\H' \cap S_{N-2} = \{q_1, \dots, q_{N-2} \}$. Let $C_i = \pi_L^{-1}(q_i)$.
Then since $S_{N-2} \cong \PP$, $C_i \sim C_j$ for all $i, j$.
Also $C_i \cap C_j = \emptyset$ if $i \ne j$. Since $\pi_L$ is projection from
$L$, $C_i \subset <L,q_i>$. Also, $C_i$ is nonsingular by Bertini's
theorem. Thus
$H \cap X = L \cup C_1 \cup \dots \cup C_{N-2}$. Thus, $\text{deg}(C_i) = b$
for all $i$. This proves \tri2.
Now by \tri2 it is clear that if $b \ge 1$ any quadric that
contains $X$ also contains $S$. Since the ideal of $S$ is generated by
quadrics \cite{ACGH, page 96}, \tri1 follows.
\tri3 follows as in the proof of Proposition 4.16.
\qed
\enddemo

\proclaim{Proposition 5.2} If $X^n \subset \P$ is a Roth variety of
degree $d = b(N-n)+1$ then
$$\CX^n = (d-b-1)^n(d-n)$$
\endproclaim

\demo{Proof}
Since $\CX = (d-n-2)H - K_X$,
$$
\align
\CX^n &= ((d-n-2)H-K_{\Xt})^n \\
&= ((d-n-2)H - (K_{\St} + \Xt))^n \cdot \Xt \\
&= ((d-n-2)H - (-(n+1)H + (N-n-2)F + bH + F)^n \cdot \Xt \\
&= ((d-b-1)H - (N-n-1)F)^n \cdot \Xt \\
&= (d-b-1)^{n-1}((d-b-1)H^n - n(N-n-1)H^{n-1}\cdot F) \cdot (bH + F) \\
&= (d-b-1)^{n-1}((d-b-1)b(N-n)+(d-b-1)-nb(N-n-1)) \\
&= (d-b-1)^{n-1}(d(d-b-1) - nb(\frac{d-1}{b} - 1) ) \\
&= (d-b-1)^n(d-n)
\qed
\endalign
$$
\enddemo

\proclaim{Corollary 5.3} If $X^n \subset \P$ is a nonsingular
nondegenerate variety of codimension $N-n \ge 2$ then $\CX$ is big
unless $X$ is the rational normal scroll $S_{1, \dots, 1}$ ($n$ 1's
with $n \ge 3$). (This is just
$\PP \times \bold P^{n-1}$ embedded into $\bold P^{2n-1}$ by the Segre
embedding).
\endproclaim

\demo{Proof}
If $\CX$ is ample then $\CX$ is
big so we can assume that $\CX$ is not ample. We can also first assume
that $X$ is linearly normal since if $\X'$ is an isomorphic projection
of $X$ then $\CX^n = \Cal C_{\X'}^n$.  Thus $X$ is a Roth variety. So
by the above proposition, $\CX^n \le 0 \iff d \le n$.  But $d =
b(N-n)+1$ and since $N-n = \sum_{i=1}^{n-1}a_i \ge n-1$ we get that
$b(n-1)+1 \le n$ and so $b=1$. Then $d \le n$ implies $(N-n)+1 \le n$
so $\sum_{i=1}^{n-1}a_i \le n-1$ and thus $a_i = 1$ for all $i$. Thus
$X$ is a Roth variety with $b=1$ contained in $S_{0,0,1,\dots,1}$ and
so the result follows when  $X$ is linearly normal. Since
the secant variety of $\PP \times \bold P^{n-1} \subset \bold P^{2n-1}$
is $\bold P^{2n-1}$, it follows that $X$ must be linearly normal.
\qed
\enddemo

\proclaim{Corollary 5.4} If $X^n \subset \P$ is a nonsingular, nondegenerate
variety of degree $d$ and codimension $N-n \ge 2$ then,
$$ H^i(X, \Cal O_X((d-n-2)H)) = 0 \qquad \text{for all \ } i \ge 1.$$
\endproclaim

\demo{Proof}
$C_X$ is  nef by Corollary 2.3 and if $C_X$ is also big then since
$K_X + \CX = (d-n-2)H$ the result follows by the Kawamata-Viehweg
vanishing theorem \cite{V}.
If $C_X$ is not big then a simple direct calculation
with the exceptional cases given in the previous corollary completes the
proof.
\qed
\enddemo

The general problem suggested here is to find conditions on $i>0$ and
$k$ in terms of the geometry of $X$ to ensure that
$H^i(X, \Cal O_X(kH)) = 0$. One context where this problem arises is in
questions involving regularity. There are two standard conjectures,
at least when $N \ge 2n+1$. They are:
\roster
\item The normality conjecture:  $X$ is $k$-normal
      if $k \ge d - (N-n)$.
\item The regularity conjecture: $X$ is $k$-regular if
      $k \ge d - (N-n) + 1$.
\endroster
The second conjecture implies the first. The regularity conjecture is
proved for (possibly singular) curves \cite{GLP} and in dimensions
2 \cite{Laz}, \cite{Pi} and 3 \cite{Ran} for nonsingular varieties. In
dimension 3, the proof needs $N \ge 9$. Much of the interest in regularity
comes from the following theorem of Mumford:

\proclaim{Theorem 5.5}
\cite{Mum1, page 99} If $X$ is $k$-regular then the homogeneous ideal
of $X$, $I(X)$ is generated by polynomials of degree $k$.
\endproclaim

By Mumford's theorem, if $X$ contains a $k+1$ secant line then
$X$ is not $k$ regular. Thus, the regularity conjecture implies
that $X$ does not have a $d+2-(N-n)$ secant line. As far as we're aware,
examples showing the sharpness of the regularity conjecture are constructed by
finding $X$ with $d+1-(N-n)$ secant lines. On the other hand, Fulton's
refined B\'ezout theorem, \cite{Fu, Example 12.3.5}
shows that $X$ cannot have a $k$ secant line if $k > d+1-(N-n)$.
(F. Zak informed us of this corollary).

Using the exact sequence:
$$
0 \to \Cal I_X(k-i) \to \Cal O_{\bold P^N}(k-i) \to \Cal O_X(k-i)
 \to 0
$$
we get that $H^{i-1}(\Cal O_X(k-i)) =
H^i(\Cal I_X(k-i))$ for all $i \ge 2$.
Thus, $X$ is $k$ regular $\iff X$ is $k-1$ normal and
$H^i(X, \Cal O_X(k-i-1))=0$ for all $i>0$. Thus the regularity
conjecture splits into two parts; the normality conjecture and the
vanishing of the higher $H^i(X,\Cal O_X(l))$'s. The second
part of the problem may seem trivial by analogy with the well
known and classical case of curves:

\proclaim{Proposition 5.6}
If $X \subset \P$ is a nondegenerate curve then
$$H^1(X, \Cal O_X(k))=0, \quad k>[(d-1)/(N-1)].$$
\endproclaim
\demo{Proof}
Use Castelnuovo's bound on the genus of a nondegenerate curve in
$\bold P^N$ to show that deg($\Cal O_X(k)) \ge 2g+1$, where $g$ is
the genus of $X$. Hence $\Cal O_X(k)$ is nonspecial.
\enddemo

Harris states \cite{HE, page 82} that an analogous vanishing is
 also true in higher dimensions i.e. that:
$$H^i(X, \Cal O_X(k))=0, \qquad k >
   \left[ \frac{d-1}{N-n}\right] , \ i>0.$$
This is false as the next examples show. First, if $Y \subset \bold P$ is
a projective variety, let $p_Y$ denote the Hilbert polynomial of $Y$.

\proclaim{Lemma 5.7}
If $A^a \subset \bold P^u$ and $B^b \subset \bold P^v$ are irreducible
$k$ normal varieties and $A \times B \subset \bold P^{uv + v + u}$
via the Segre embedding then
\roster
\item $A \times B$ is $k$ normal in $\bold P^{uv + v + u}$.
\item $p_{A \times B}(k) = p_A(k)p_B(k)$ and so
deg($A \times B) = \binom{a+b}{b} \text{deg}(A) \text{deg}(B)$.
\endroster
\endproclaim
\demo{Proof}
We have the commutative diagram
$$
\CD
H^0(\bold P^{uv+v+u}, \Cal O_{\bold P^{uv+v+u}}(k)) @>\a>>
 H^0(A \times B, \Cal O_{A \times B}(k)) \\
@V{\gamma}VV @| \\
H^0(\bold P^u \times \bold P^v, \Cal O_{\bold P^u \times \bold P^v}(k))
@>\b>> H^0(A \times B, \Cal O_{A \times B}(k))
\endCD
$$
It is easy to see that $\gamma$ is surjective. Hence $A \times B$ is
$k$ normal $\iff \b$ is surjective.
If $H^0(\bold P^u, \Cal O_{\bold P^u}(k)) @>\b_1>> H^0(A, \Cal O_A(k))$
and $H^0(\bold P^v, \Cal O_{\bold P^v}(k)) @>\b_2>> H^0(B, \Cal O_B(k))$ are
the natural restriction maps then by the Kunneth formula we obtain
a commutative diagram where the vertical arrows are isomorphisms:
$$
\CD
H^0(\bold P^u \times \bold P^v, \Cal O_{\bold P^u \times \bold P^v}(k))
@>\b>> H^0(A \times B, \Cal O_{A \times B}(k)) \\
@VVV @VVV \\
H^0(\bold P^u, \Cal O_{\bold P^u}(k)) \otimes
H^0(\bold P^v, \Cal O_{\bold P^v}(k)) @>{\b_1 \otimes \b_2}>>
H^0(A, \Cal O_A(k)) \otimes H^0(B, \Cal O_B(k))
\endCD
$$
And thus if $\b_1$ and $\b_2$ are surjective, then so is $\a$. \tri1 follows.

Using the Kunneth formula, for $k \gg 0$
$$
\align
p_{A \times B}(k) &= h^0(A \times B, \Cal O_{A \times B}(k)) \\
  &= h^0(A, \Cal O_A(k))h^0(B, \Cal O_B(k)) = p_A(k)p_B(k).
\endalign
$$
This proves \tri2.
\qed
\enddemo

\subheading{Example 5.8}
Let $A \subset \bold P^2$ be a nonsingular plane curve of
degree $d_A$. Let $X^n = A \times \bold P^{n-1} \subset \bold P^{3n-1}$
via the Segre embedding and assume $n \ge 2$. By the above lemma, $X$ is
projectively normal and the degree of $X$, $d = nd_A$. By the
Kunneth formula,
$$
H^1(X,\Cal O_X(kH)) \cong H^1(A, \Cal O_A(kH)) \otimes H^0(
\bold P^{n-1}, \Cal O_{\bold P^{n-1}}(kH)).
$$
Since $H^1(A, \Cal O_A(kH)) \cong H^2(\bold P^2, \Cal O_{\bold P^2}(k-d_A))$,
if $0 \le k \le d_A-3$ it follows that  $H^1(X, \Cal O_X(kH)) \ne 0$.
Harris' claim implies that
$$
H^1(X, \Cal O_X(kH)) = 0, \qquad k \ge \left[ \frac{nd_A-1}{2n-1}\right]
$$
But since $d_A-3 > \left[ \frac{nd_A-1}{2n-1}\right]$ for $d_A$ sufficiently
large, $X$ provides a counterexample.

Note also that by Proposition 3.11 \tri5, $X$ is projectively normal even
though a general hyperplane section of $X$ is not.

We now want to determine the normal bundle of the line in a Roth
variety.
\proclaim{Proposition 5.9}
Let $X^n \subset \P$ be a Roth variety
of degree $d=b(N-n)+1$ and $S$ a rational normal scroll
$S_{0,0.a_1, \dots, a_{n-1}}$ with all $a_i \ge 1$, $\sum_{i=1}^{n-1}a_i
=N-n$, and singular locus $L$ such that $L \subset X \subset S$. Then
$N_{L/X} = \Cal O_{L}(1-ba_1) \oplus \dots \oplus \Cal O_{L}(1-ba_{n-1})$.
In particular, $c_1(N_{L/X}) = n-d$.
\endproclaim

\demo{Proof}
First we'll prove the formula for the Chern class $c_1(N_{L/X})$ by
induction on $n$.
Suppose that $n=2$. Then $c_1(N_{L/X}) = L^2$. Using the
notation of Theorem 3.8 and the intersection theory developed in that
section we see that
$$
\align
L^2 &= \pi_2|_{\Xt}^{-1}(L)^2  = (\PP \times L)^2 \cdot \Xt \\
    &= (H - (N-2)F)^2 \cdot (bH + F) = bH^3 + (1- 2b(N-2))H^2 \cdot F \\
    &= 1 - b(N-2) = 2- d.
\endalign
$$

Now suppose $n \ge 3$. Let $H$ be a generic hyperplane containing $L$ and
let $Y = X \cap H$. Then by Proposition 5.1 \tri3, $Y$ is a Roth variety.
We have the standard exact sequence
$$
0 \to N_{L/Y} \to N_{L/X} \to N_{Y/X}|_{L} \to 0.
$$
Since $N_{Y/X}|_{L} = \Cal O_L(1)$ and by the induction hypothesis
$c_1(N_{L/Y}) = n-1-d$ it follows that $c_1(N_{L/X}) = n - d$.

Using the notation of Proposition 4.11 (and the paragraph preceeding it)
we now want to study the diagram
$$
\CD
E @>f_2|_E>> Z \\
@Vf_1|_EVV  @. \\
L @.
\endCD
$$
where $E = \bold P(N_{L/X})$, $Z = S_{a_1, \dots a_{n-1}}$ and
$f_2|_E$ is surjective. To simplify notation, let $f = f_2|_E$
and $g= f_1|_E$. For each $p \in E$, $f(g^{-1}(p))$ is an
$n-1$ dimensional linear space on the scroll $Z$. Thus we can write
$Z$ as a $\bold P^{n-2}$ bundle $g_1:Z \to \PP$ and
obtain a commutative diagram
$$
\CD
E @>f>> Z \\
@VgVV @Vg_1VV \\
L @>f_1>> \PP
\endCD
$$
as in the last paragraph of the proof of Claim 3.6.
In particular, $f$ is finite and has the same degree as $f_1$.
Using Proposition 5.1 \tri2, it follows from the $n=2$ case that
the degree of $f_1$ is $b$.
If $\Cal E = \Cal O_{\PP}(a_1) \oplus \dots \oplus \Cal O_{\PP}(a_{n-1})$ then
$Z = \bold P(\Cal E^*)$.

Now by \cite{Ha, Chapter II, Proposition 7.12}, to give a morphism
$E \to Z$  over $\PP$ as above is equivalent to giving a line bundle
$\Cal L$ on $E$ and a surjection
$(g \circ f_1)^*\Cal E \to \Cal L$.
Thus we obtain a surjection $f_1^* \Cal E \to g_*\Cal L$. Since $f$ is
surjective, this surjection is in fact an isomorphism.
Thus $g_*(\Cal L) = N_{L/X}^* \otimes \Cal O_L(k)$ for some integer
$k$. Also, since the degree of $f_1$ is $b$,
$f_1^* \Cal E = \Cal O_L(ba_1) \oplus \dots \oplus \Cal O_L(ba_{n-1})$.
Thus, if $N_{L/X} = \Cal O_L(b_1) \oplus \dots \Cal O_L(b_{n-1})$ then
we can reorder the $b_i$ to conclude that $b_i = k - ba_i$ for all $i$.
Since $\sum_{i=1}^{n-1}b_i = n-d$, $\sum_{i=1}^{n-1}a_i = N-n$ and
$d = b(N-n) + 1$ we conclude that $k = 1$ which gives the required
result.
\qed
\enddemo

\proclaim{Proposition 5.10} If $S$ is a rational
 normal scroll
$S_{a_0, \dots, a_{n-1}} \subset \P$ with $\sum a_i = N-n+1$
and all $a_i \ge 0$ then for $S_{a_0} \subset S$,
$$N_{S_{a_0}/S} = \Cal O_{S_{a_0}}(a_0 - a_1) \oplus \dots \oplus
\Cal O_{S_{a_0}}(a_0-a_{n-1}).$$
\endproclaim

To prove this we use the following fact \cite{Fu, Appendix B.7.4}:
\proclaim{Fact 5.11}
If $Z$ is a nonsingular variety and $Y_i \subset X$ for $1 \le i \le r$
are nonsingular subvarieties of codimension $d_i$ such that $X = \cap Y_i$ is
a nonsingular variety of codimension $\sum d_i$ then
$N_{X/Z} = \oplus_{i=1}^{r}(N_{Y_i/X})|_{X}$.
\endproclaim

\demo{Proof of Proposition}
We will apply the above fact with $Z = S$, $r=n-1$, and
$Y_i = S_{a_0, a_1, \dots, \hat a_i, \dots, a_{n-1}}$.
Then $X = \cap Y_i = S_{a_0}$.
Let
$\Cal E =  \Cal O_{\PP}(a_0) \oplus \dots
\oplus \Cal O_{\PP}(a_{n-1})$ and let
$\pi_1:\bold P(\Cal E^*) \to \PP$ be a $\bold P^{n-1}$ bundle with diagram:
$$
\CD
\bold P(\Cal E^*) @>\pi_2>> S \subset \P \\
@V\pi_1VV @. \\
\PP @.
\endCD
$$
such that the morphism $\pi_2$ is given by the
complete linear series $|\Cal O_{\bold P(\Cal E^*)}(1)|$. Let H be the
pullback of a hyperplane in $\P$ via $\pi_2$ and let $F$ be a fiber of $\pi_1$.
Then $H^n = \sum_{j=0}^{n-1}a_j$, $H^{n-1}\cdot F = 1$ and $F^2 = 0$.
Since $Y_i \cdot F \cdot H^{n-2} = 1$ and
$Y_i \cdot H^{n-1} = -a_i + \sum_{j=0}^{n-1}a_i$ it follows that
$Y_i \sim H - a_iF$. Thus  $N_{Y_i/S} = \Cal O_{Y_i}(H-a_iF)$ and
so $N_{Y_i/S}|_{X}= \Cal O_{X}((H-a_iF) \cdot X) =
\Cal O_{X}((H-a_iF) \cdot (H-a_1F) \cdot \dots \cdot (H-a_{n-1}F)) =
\Cal O_{X}(a_0-a_i)$. The result now follows by the above fact.
\qed
\enddemo

Finally, to determine the associated scroll of the
Roth variety $Y$ in Proposition 5.1 \tri3 we need to determine
the generic hyperplane section of a rational normal scroll
$S_{a_1, \dots, a_n}$. We first heard the statement of the next
proposition from Miles Reid but we could not find a proof in the
literature so it is given here.
\proclaim{Proposition 5.12}
\roster
\item
There is a surjection
$$
\Cal O_{\PP}(a_1) \oplus \dots \oplus \Cal O_{\PP}(a_n) \to
\Cal O_{\PP}(b_1) \oplus \dots \oplus \Cal O_{\PP}(b_m)
$$
$\iff m \le n$ and after reordering so that
$a_1 \le \dots \le a_n$ and $b_1 \le \dots b_m$ then
$b_i \ge a_i$ and if $(a_1, \dots, a_i) \ne (b_1, \dots, b_i)$ then
also $b_i \ge a_{i+1}$.
\item Suppose all $a_i, b_i \ge 1$. Then
$S_{b_1, \dots, b_{n-1}}$ is a
hyperplane section of $S_{a_1, \dots, a_n}$
$\iff \sum a_i = \sum b_i$ and there is a surjection
$$
\Cal O_{\PP}(a_1) \oplus \dots \oplus \Cal O_{\PP}(a_n) \to
\Cal O_{\PP}(b_1) \oplus \dots \oplus \Cal O_{\PP}(b_{n-1}).
$$
\endroster
\endproclaim

\demo{Proof}
We first prove the ``$\implies$'' direction for \tri1.
A surjection as above is given by an $m \times n$ matrix $T$ whose
$i,j^{th}$ entry is a polynomial  $T_{i,j} \in H^0(\PP,\Cal O_{\PP}(b_i-a_j))$
such that evaluation at each point $p \in \PP$ gives a rank $m$ matrix
$T(p)$. Suppose $b_k < a_k$ for some $k$.  Then $T_{i,j} = 0$ for
all $i \le k$ and $j \ge k$. It follows that $M(p)$ has rank $\le m-1$
for all $p \in \PP$. This contradiction shows that $b_k \ge a_k$.
Suppose that $(a_1, \dots, a_k) \ne (b_1, \dots, b_k)$ and
$b_k < a_{k+1}$ for some $k$. Then $T_{i,j}=0$ for
all $i \le k $ and $j \ge k+1$. Since $(a_1, \dots, a_k) \ne (b_1, \dots, b_k)$
, the determinant of the first $k \times k$ minor of $T$ is a nonconstant
polynomial or $0$. Hence it vanishes for some $p$ in $\PP$. For such a
value of $p$, $T(p)$ has rank $\le m$. This completes the proof.

To prove the other direction of \tri1, consider the matrix $T$
whose nonzero entries are: $T_{i,i} := x_0^{b_i-a_i}$ and
$T_{i,i+1} := x_1^{b_i - a_{i+1}}$ if $b_i \ge a_{i+1}$ and $0$
otherwise. It is easy to see that $T(p)$ has rank $m$ for all $p$ and
hence determines a surjective map as required.

Let $\Cal E = \Cal O_{\PP}(a_1) \oplus \dots \oplus \Cal O_{\PP}(a_n)$
and $\Cal F = \Cal O_{\PP}(b_1) \oplus \dots \oplus \Cal O_{\PP}(b_{n-1})$.
Let
$\pi_1:\bold P(\Cal E^*) \to \PP$ be a $\bold P^{n-1}$ bundle with diagram:
$$
\CD
\bold P(\Cal E^*) @>\pi_2>> S_{a_1,\dots,a_n} \subset \P \\
@V\pi_1VV @. \\
\PP @.
\endCD
$$
with the morphism $\pi_2$ is given by the
complete linear series $|\Cal O_{\bold P(\Cal E^*)}(1)|$. Let H be the
pullback of a hyperplane in $\P$ via $\pi_2$ and let $F$ be a fiber of $\pi_1$.

Given a surjection $\Cal E \to \Cal F$, we obtain an inclusion
$\bold P(\Cal F^*) \subset \bold P(\Cal E^*)$ over $\PP$ such that
$\Cal O_{\bold P(\Cal E^*)}(1)|_{\bold P(\Cal F^*)} =
\Cal O_{\bold P(\Cal F^*)}(1)$. Thus $H^{n-1} \cdot  \bold P(\Cal F^*)
= (H|_{\bold P(\Cal F^*)})^{n-1} = \sum b_i$. Also $H^{n-2} \cdot F
\cdot \bold P(\Cal F^*) = 1$. Thus $\bold P(\Cal F^*) \sim H$ and so
$\pi_2(\bold P(\Cal F^*)) = S_{b_1, \dots, b_{n-1}}$, a hyperplane
section of $S_{a_1, \dots, a_n}$.

Conversely, given that $S_{b_1, \dots, b_{n-1}}$ is a hyperplane
section of $S_{a_1, \dots, a_n}$ we obtain an inclusion
$\bold P(\Cal F^*) \subset \bold P(\Cal E^*)$ over $\PP$
such that $\Cal O_{\bold P(\Cal E^*)}(1)|_{\bold P(\Cal F^*)} =
\Cal O_{\bold P(\Cal F^*)}(1)$.
By \cite{Hart, Chapter II, Proposition 7.12} this gives a surjection
$(\pi_1|_{\bold P(\Cal F)})^* \Cal E \to \Cal O_{\bold P(\Cal F^*)}(1) $.
 Pushing this forward to $\PP$
we obtain a surjection $\Cal E \to \Cal F$.
\qed
\enddemo

\proclaim{Theorem 5.13} \cite{Ha1, Section 3}
$S_{b_1, \dots, b_n}$ is a degeneration of $S_{a_1, \dots, b_n}$ in
the Chow variety of degree $N-n+1$ $n$-folds in $\P$ $\iff$
after reordering so that $a_1 \le \dots \le a_n$ and $b_1 \le \dots \le b_n$
then $\sum_{i=1}^{\a}b_i \le \sum_{i=1}^{\a}a_i$ for all $\a$.
\endproclaim

Using the above two results one can determine the generic hyperplane
section of a given rational normal scroll
e.g. the generic  hyperplane section of $S_{5,9,11,15}$ is $S_{12,13,15}$.

Finally, we have some concluding remarks.
\subheading{Remark 5.14} Analogous to the study of Roth varieties, we can
consider $X^n \subset \P$ a nonsingular, nondegenerate variety
such that there is an $n+1$ dimensional rational normal scroll
$S= S_{0,a_1,\dots,a_n}$ with all $a_i \ge 1$, $\sum a_i = N-n$ and singular
point $p$ such that $p \subset X \subset S$. {\sl Mutatis mutandis},
all the results of Section 3 for Roth Varieties (except those
expressing properties of the line $L$ in the Roth variety) hold for the above
varieties and can be proved using the techniques of that section.
In particular, Theorem 3.8 and Theorem 3.18 hold.

\subheading{Remark 5.15} Let $X^n \subset \P$ be a nonsingular, nondegenerate
variety of codimension $N-n \ge 2$.
The most obvious question left unresolved in this work
is to understand when the linear system determined by the double-point
divisor on $X$, $|\CX|$ is very ample. If $X$ is not a Roth variety
and $|\CX|$ is not very ample then $|\CX|$ separates points but does not
separate tangent directions. We don't know if such examples exist.
(A possible candidate for an example would be a variety as in the preceding
remark where $|\CX|$ would not separate tangent directions at the
singular point $p$ of the associated rational normal scroll.)
An even more interesting question is to understand when
 tangent directions can be separated using only geometrically constructed
sections of $|\CX|$ i.e. sections of the form $\CXL$ for some center
of projection $\L$.

\bigskip
\bigskip
\noindent
Author's address:
Math Dept., Columbia University, New York, N.Y., 10027. \break
E-mail: ilic\@math.columbia.edu

\end